\documentclass[a4paper,10pt]{article}

\def\q{\quad}

\def\cut{\hfill\break}  \def\newline{\cut}
\def\h#1{\hbox{#1}}

 \def\eps{\epsilon}
\def\be{\beta}
\def\a{\alpha}  
  
\def\vp{\varphi}

\def\calC{{\h{\cal C}}} \def\C{{\calC}}

\def\II{{Id}}
 
\def\RR{{\mathbb R}} 
\def\ZZ{{\mathbb Z}} 
\def\NN{{\mathbb N}}

\def\diam{\hbox{\rm diam}}

   \def\Ric{\h{\rm Ric}}

\def\osc{\hspace{0.02cm}\textrm{osc}_M\hspace{0.02cm}}

\def\V{\frac1V}   

\def\gij{g_{i\bar j}}

\def\oKE{\o_{KE}}

\def\AutMJ{\h{\rm Aut}(M,J)}

\def\GL2nR{\h{$GL(2n,\RR)$}}
\def\Sp2nR{\h{$Sp(2n,\RR)$}}

\def\dbz{d\bar z}

\def\dzidzjb{dz^i\w\dbz^j}

\def\ddt{\frac d{dt}}

\def\all{\forall}

\def\precpt
{\hbox{$\mskip3mu\mathhalf\subset\raise0.92pt\hbox{$\mskip-10mu\!\!\!\!
\mathfootnote\subset$}\mskip5mu$}}

\def\supsetnoteq{\hbox{$\mskip3mu\supset\raise-5.97pt
\hbox{$\mskip-10mu\!\!\!\scriptstyle\not=$}\mskip8mu$}}

\def\subsetnoteq{\hbox{$\mskip3mu\subset\raise-5.97pt
\hbox{$\mskip-11mu\!\!\scriptstyle\not=$}\mskip8mu$}}

\def\sseq{\subseteq}

\def\*{\star}

\def\w{\wedge}

\def\D{\Delta}
\def\N{\nabla} 
\def\dbar{\bar\partial}
\def\del{\partial} 
\def\ddbar{\partial\dbar}

\def\intM{\int_M}  \def\intm{\int_M}
\def\V{\frac1V}
\def\VintM{\V\int_M}

 \def\ha{\frac12}      

\def\ra{\rightarrow}

\def\MA{Monge-Amp\`ere }
   \def\MAes{Monge-Amp\`ere equations}

\def\K{K\"ahler }

\def\KE{K\"ahler-Einstein }  \def\KEno{K\"ahler-Einstein}

\def\CtwoM{\C^2(M)}

\def\Ctwobe{\C^{2,\beta}}

\def\CtwobeMo{\C^{2,\beta}(M,\omega)}

\def\Cinf{\h{\cal C}^\infty}

\def\CinfM{\Cinf(M)}

\def\LpMo{L^p(M,\o)}
\def\LinfMo{L^\infty(M,\o)}

\def\n#1#2{||#1||_{#2}}

\def\on{\o^n}

\def\ovp{\o_{\vp}}

\def\vpt{\vp_{\!t}}

\def\ovpt{\o_{\vpt}}

\def\ovpn{\o_{\vp}^n}
\def\ovptn{\o_{\vpt}^n}

\def\Ho{{\call H_{\o}}}

\def\ric#1{\hbox{\rm Ric}^{(#1)}}
\def\ricm#1{\hbox{\rm Ric}^{(-#1)}}
\def\Ric{\hspace{0.025cm}\textrm{Ric}\hspace{0.025cm}}

\def\Hc{K\!a_{c_1}}

\def\hc#1{K\!a_{c_1}^{(#1)}}

\def\hcbar#1{\overline{K\!a_{c_1}^{(#1)}}}

\def\HcG{K\!a_{c_1}(G)}

\def\Hmuc{K\!a_{\mu c_1}}

\def\O{\Omega}
\def\rico#1{\hbox{\rm Ric}_{\Omega}^{(#1)}}

\def\RicmO{\hbox{\rm Ric}^{(-1)}_\O}

\def\Ho{K\!a_{\Omega}}

\def\ho#1{K\!a_{\Omega}^{(#1)}}

\def\HoG{K\!a_{\Omega}(G)}

\def\Aub{\h{\rm Ric}}
\def\Aubeps{\h{\rm Ric}_\epsilon}
\def\Aubepsl{\h{\rm Ric}_\epsilon^{(l)}}
\def\Aubepslj{\h{\rm Ric}_\epsilon^{(l_j)}}
\def\Aubepsj#1{\h{\rm Ric}_\epsilon^{(#1)}}

\def\Fmu{F_\mu}

\def\MA{Monge-Amp\`ere }

\def\oone{\o_{\vp_1}}
\def\otwo{\o_{\vp_2}}

\def\Phil{\Phi_l}

\def\vpone{\vp_1}
\def\vptwo{\vp_2}

\def\ovpone{\o_{\vpone}}
\def\ovptwo{\o_{\vptwo}}
\def\Dvpt{\D_{\ovpt}}
  \def\V{\frac1V}

\def\fovp{f_{\ovp}}

\def\ol{\o_l}

\def\oln{\o_l^n}

\def\ovpln{\o_{\Sigma^l_{j=1}\vp_j}^n}

\def\Phil{\Phi_l} 

\def\Fmu{F_\mu}

\def \Mark #1{\xdef #1{(\the \eqncount)}#1}

\usepackage{graphicx}
\usepackage{amssymb}
\usepackage{amsfonts}
\usepackage{amsmath}
\usepackage{amsthm}

\usepackage[colorlinks=true]{hyperref}
\def\o{\omega}
\def\O{\Omega}
\def\i{\sqrt{-1}}

\newtheorem{remark}{Remark}[section]
\newtheorem{lemma}{Lemma}[section]
\newtheorem{definition}{Definition}[section]

\newtheorem{notation}{Notation}[section]
\newtheorem{proposition}{Proposition}[section]
\newtheorem{corollary}{Corollary}[section]
\newtheorem*{claim}{Claim}
\newtheorem{theorem}{Theorem}
\newtheorem{conjecture}{Conjecture}
\newtheorem{theorem-others}{Theorem}[section]
\title{Ricci iterations on K\"ahler classes}
\author{Julien Keller}

\begin{document}

\maketitle

\begin{abstract}
In this paper we consider the dynamical system involved by the Ricci operator
on the space of \K metrics. A.\ Nadel has defined an iteration scheme given by the Ricci operator for Fano manifold and asked whether it has some nontrivial periodic points. 
First, we prove that no such periodic points can exist. We define the inverse of the Ricci operator and consider the dynamical behaviour of its iterates for a Fano \KE manifold. In particular we show that the iterates do converge to the K\"ahler-Ricci soliton for toric manifolds. Finally, we define a finite dimensional procedure to give an approximation of \KE metrics using this iterative procedure and apply it for $\mathbb{P}^2$ blown up in 3 points.
\end{abstract}

\section*{Introduction}

\smallskip

Let $M$ be a compact \K manifold of complex dimension $n$. For any \K metric 
$g$ we
let $\o:=\o_g=\sqrt{-1}/{2\pi}\cdot\gij(z)\dzidzjb$ denote its corresponding \K form,
a closed positive (1,1)-form on $M$. 
For a \K form $\o$ such that $[\o]=c_1(M)$ consider the space of strictly  $\o$-plurisubharmonic potentials
$$
K\!a_{[\o]}=\{\vp\in\CinfM\,:\, \o+\sqrt{-1}\ddbar\vp >0\},
$$
and the space $\Hc$ of \K forms cohomologous to $c_1(M)$.

In \cite{Na}, Nadel considers for a Fano manifold $M$, iterations on $\Hc$
defined
inductively using the Ricci operator as follows. Let $\hc0=\Hc$ denote
the set of all smooth forms in $c_1(M)>0$ and let
$\ric0$ be the identity operator. For any \K
metric we let $$\Ric(\omega)=-\ddbar\log\det(\gij)$$ denote the
Ricci form of $\o$. It is well-defined globally and lies again in
the $c_1(M)$ class. If it is positive we let $\ric2(\omega)$ denote its
Ricci form, and in a similar fashion we define higher powers of
the operator as long as the positivity is preserved. The motivation for this construction comes from the simple fact that, when they exist, \KE
metrics are by definition fixed points for this iteration process.
Nadel asked whether these are all periodic points and proved the absence of 
periodic points of order 2 and 3.
Furthermore he raised the question whether the existence of \KE metrics could be
related to this iteration procedure. This question is also very natural. Actually, as we 
will explain later, one can define the operator $\Ric^{(-1)}$ using Calabi-Yau theorem and see its iterations as a kind of naive {\it discretization} of the (normalised) K\"ahler-Ricci flow
\begin{equation}\label{ricci-flow}
\frac{\partial w_t}{\partial t} = -\Ric(\omega_t) + \omega_t
\end{equation}
Thus, we expect some similar properties for the K\"ahler-Ricci flow and the iterations of the operator $\Ric^{(-1)}$.  \bigskip

\par We explain now the organization of this paper. We will show that some natural energy functionals are decreasing along these iterations. This will give us a simple proof of the non existence of (non trivial) periodic points and thus answer Nadel's first question. Then we study the question of the existence of periodic points of infinite order. The behaviour of our dynamical system is closely tied with the existence of \KE metrics on the Fano manifold. We generalize Nadel's iteration scheme and define a family of natural operators $\Aubeps, \,\eps\ge0$ and see that their behaviour is particularly simple when the manifold is Einstein and $\eps<1$. In Section \ref{toric}, we investigate the case of Fano toric manifolds for which we can prove that the iterates of $\Ric^{(-1)}$ do always converge towards a K\"ahler-Ricci soliton. This also gives a new proof of the existence of K\"ahler-Ricci solitons on Fano toric manifolds. We shall emphasise the fact that we don't use any flow or continuity method, contrarily to \cite{Zh1,WZ}. Then we relate the iteration procedure to the notion of canonically balanced metric studied by Donaldson in \cite{D3}. This gives us an approximation procedure in a finite dimensional setup of the \KE metric when it does exist a priori on a Fano manifold. We apply our techniques to the case of $\mathbb{P}^2$ blown up in $3$ points and give a numerical approximation of the \KE metric living on it. Finally we discuss the case of the manifolds with negative first Chern class and also non canonical classes.

\section{Positive and negative Ricci iterations}\label{SectionRicci}

Our first observation lies in the fact that Nadel's
construction can actually be reversed. Let $\alpha$ be any form
representing the $c_1(M)$ class. By the Calabi-Yau theorem
\cite{Ya} there exists a unique \K form in $\Hc$, which we denote by
$\omega_1$, whose Ricci form equals $\alpha$. We define the {\it
inverse Ricci operator} by $$\ricm1(\a):=\o_1.$$ Similarly we define
higher powers $\ricm k:=\ricm1\circ\cdots\circ\ricm1$, making
repeated use of the Calabi-Yau theorem.  This leads us to consider for a given
K\"ahler form $\omega_0$ the sequence
\begin{equation}\label{iter_einstein}
\Ric(\omega_{j+1}) = \omega_j 
\end{equation}
for an integer $j\geq 1$. This can be also written
$$ \frac{\omega_{j+1} -\omega_{j}}{(j+1)-j} = - \Ric(\omega_{j+1}) + \omega_{j+1}$$
which justifies naively why our iterative procedure (\ref{iter_einstein}) should have some similarities
with the K\"ahler-Ricci flow (\ref{ricci-flow}).
\par For a positive $k$ we will denote
$\hc k$ the maximal domain of definition of $\ric k$. We
obtain naturally a filtration of $\hc0$. We also let $\hcbar k$ denote
the set of all metrics in $\hc {k-1}$ whose image under $\ric {k-1}$ is nonnegative.

The advantage of this
construction lies in the fact that iterations {\it improve} the
positivity of $\alpha$. Indeed, if we consider the filtration
$\{\hc j\}_{j\ge0}$ to measure positivity in some sense, then the
image of $\ricm k$ lies in $\hc k$ and we may iterate the inverse
Ricci operator to any desired power. Note that this yields
homeomorphisms
$$\ric k:(\hc l,||\cdot||_{\C^{m,\beta}}) \ra(\hc
{l+k},||\cdot||_{\C^{m-2k,\beta}})$$ for all $\beta\in(0,1)$,
$m,l\in\NN\cup\{0\},k\in\ZZ$ such that $l+k\ge0$ and $m\ge2k$.
We let $\hc{\infty}$ denote the set of all $L^\infty(M)$ limits
$\lim_{k\ra\infty}\ricm k\alpha$, for $\alpha\in\hc0$ (when they
exist). 

Finally, let $G$ be any connected compact subgroup 
of the group $\AutMJ$ of holomorphic diffeomorphisms of $(M,J)$. Denote by
$\HoG$ the space of $G$-invariant \K forms in $\Ho$. Such forms exists as can be seen
by averaging over orbits of $G$ with respect to the Haar measure of $G$. We remark that $\ricm1$ maps $\HcG$ into itself.

\section{Energy functionals on the space of K\"ahler potentials}
\label{SectionEnergyFunctionals}
In the Sections \ref{F} and \ref{K}, we present some properties of some well-known energy functionals. Let $\O\in H^2(M,\RR)$ denote a \K class. 
We call a function $A:\Ho\times\Ho\ra\RR$ an {\it energy functional} if it 
is zero on the diagonal. By an {\it exact energy functional} we will mean one
which satisfies in addition the cocycle condition (see \cite{Ma})
$$A(\o_1,\o_2)+A(\o_2,\o_3)=A(\o_1,\o_3).$$

In all the paper, $V$ will denote the volume of the manifold with respect to $[\omega]$, i.e $V=\int_M \frac{\o^n}{n!}$.
The energy functionals
$I,J$, introduced by Aubin in \cite{A2}, are defined for each pair 
$(\o,\o_{\vp}:=\o+\i\ddbar\vp)$ by
\begin{equation}
I(\o,\ovp)=\frac{1}{V}\int_M\i\del\vp\w\dbar\vp\w\sum_{i=0}^{n-1}\o^{i}\w\ovp^{n-1-i}
=\frac{1}{V}\int_M\vp(\on-\ovpn)\label{Ieq},
\end{equation}
\begin{equation}
J(\o,\ovp)=\frac{1}{V(n+1)}\int_M\i\del\vp\w\dbar\vp\w\sum_{i=0}^{n-1}(n-i)\o^{i}\w\ovp^{n-1-i}.\label{Jeq}
\end{equation}
We note some of their basic properties for which we refer the reader to 
\cite{A2,Si,T2}.
Note that $I,J$ and $I-J$ are all nonnegative and equivalent. One may also define
them via a variational formula. 
Connect each pair
$(\o,\o_{\vp_1}:=\o+\i\ddbar\vp_1)$ with a piecewise smooth path
$\{\ovpt\}$. Then, e.g. for $I-J$, we have for any such path
\begin{equation}
(I-J)(\o,\o_{\vp_1})=-\frac{1}{V}\int_{[0,1]\times M}\vpt\D_t\dot\vpt\ovptn\w dt.\label{IminusJeq}
\end{equation}

\subsection{The $F_\mu$ functional}\label{F}

Let 
\begin{equation}
F^0(\omega_0,\vp)=-(I-J)(\omega_0,\ovp)-\V\intM\vp\ovp^n.\label{Fzeroeq}
\end{equation}
For a \K manifold of positive or negative Chern class define the following functional
on $\Hmuc\times\Hmuc$
$$
F_\mu(\omega,\omega_\vp)=F^0(\o,\vp)-\mu\log\left( \frac{1}{V}\intM e^{h_{\o}-\mu\vp}{\o}^n \right),
$$
with $\mu=\pm1$ respectively.
For a \K manifold of zero Chern class define on $\Ho\times\Ho$ the functional
$$
F_0(\omega,\ovp)=F^0(\o,\vp)+\frac{1}{V}\intM \vp e^h\o^n.
$$
The critical points of the functionals $\Fmu$ are the \KE metrics. They are absolute minima \cite{Di,T2}. Indeed the second variation of $F$ at a critical point in the direction of the plane
spanned by $\psi_1,\psi_2\in T_{\vpt}Ka_{[\omega]}$ is given by
\begin{equation}
\V\intM[\ha g_{\vpt}(\N\psi_1,\N\psi_2)-\mu\psi_1\psi_2]\ovptn,\q \mu=\pm1,0.\label{FSecondVariationeq}
\end{equation}

This is seen to be a strictly positive $(0,2)$-tensor on $Ka_{[\omega]}$ for $\mu\le-1$. For $\mu=1$
it is nonnegative and vanishes precisely when $\psi_1$ and $\psi_2$ are proportional 
and eigenfunctions of eigenvalue $-1$ of $\D_{\dbar}$ (see \cite[p.64]{T2}). In that case this infinitesimal variation corresponds to holomorphic automorphisms and to moving within
the set of \KE forms.

\subsection{K-energy and $E_k$ functionals}\label{K}
The Chen-Tian energy functionals $E_k,\,\,k=0,\ldots,n$, are
defined in a similar manner by  
\begin{eqnarray}
E_k(\omega,\omega_{\vp_1}) &=& \frac{(k+1)}{V}
\int_{[0,1]\times M}
\D_{\vp_t}\dot\vp_t\Ric(\o_{\vp_t})^k\w\ovpt^{n-k}\w dt \nonumber \\
&&-\frac{(n-k)}{V} \int_{[0,1]\times M}
\dot\vp_t(\Ric(\o_{\vp_t})^{k+1}-\mu_k\ovpt^{k+1})\ovpt^{n-1-k}
dt, \\
\mu_k&=&\frac{c_1(M)^{k+1}\cup[\o]^{n-k-1}([M])}{[\o]^n([M])}.
\end{eqnarray}
This gives rise to
well-defined {\it exact} energy functionals independent of the choice of path \cite{CT}.
The {\it K-energy}, $E_0$, was introduced by Mabuchi \cite{Ma}.
The following formula is taken from \cite[Section 7.2]{T0}. For the derivation (for any \K class) we refer to \cite{Ch} where an equivalent expression is given.
\begin{proposition} 
Let $h$ be a function satisfying ${\rm Ric}\o-\mu\o=\i\ddbar h$. One has 
$$
E_0(\o,\ovp)=
\V\intM\log\frac{\ovpn}{\on}\ovpn-\mu(I-J)(\o,\ovp)+\V\intM h(\on-\ovpn).
$$
\end{proposition}

Remark that the $E_k$ functionals vanish
on pairs joined by a one parameter subgroup of automorphisms through the identity
\cite[Corollary 5.5]{CT}.

\begin{definition}
We say that an exact functional $T$ is bounded from below if $T(\o,\ovp)\ge C$ for every
$\ovp\in\Ho$. We say it is proper on $\HoG$ (in the sense of Tian) if there exists a function $\rho:\RR\ra\RR$ 
satisfying $\lim_{s\ra\infty}\rho(s)=\infty$ such that $T(\o,\ovp)\ge\rho((I-J)(\o,\ovp))$
for every $\ovp\in\HoG$. 
\end{definition}

This is well-defined, in other words
depends only on $[\o]$ since the failure of $I-J$ to satisfy the cocycle condition is 
under control with respect to the two base metrics,
$$
(I-J)(\o,\otwo)-(I-J)(\oone,\otwo)=
(I-J)(\o,\oone)-\V\intM\vp_1(\otwo^n-\oone^n).
$$

\subsection{A lower bound for the energy functionals}
We now study the iterations of the $\Ric^{-1}$ operator for Fano manifolds by analysing
the behaviour of the functionals that we have introduced. Firstly, we recall a result of non-negativity of Bando and Mabuchi,

\begin{theorem-others} \label{BMLowerboundThm}
\cite[Theorem A]{BM},\cite[Theorem 1]{Ba},\cite[Theorem 1.2]{SW}. 
Let $(M,\oKE)$ be a Fano \KE manifold. 
Then for $i=0,1$,
$$E_i(\o_{KE},\o)\ge0$$ for all $\o\in\Hc$ with equality if and only
if $\o$ is \KEno. In that case there exists a holomorphic automorphism
homotopic to the identity $h$ such that $h^\star\oKE=\o$.
\end{theorem-others}

We give a sketch of the proof for $i=0$ with an emphasis on 
the features that will be useful in later sections.
Consider the deformation $\{\ovpt\}\sseq\Hc$ constructed
from two paths, solutions of the following \MA equations
\begin{eqnarray}
\ovpt^n&=& {e^{tf+c_t}\o^n},t\in[0,1] \\
       &=&  e^{f-(t-1)\vpt}\o^n,t\in[1,2]\label{MApathseq}
\end{eqnarray}
where $\Ric\o-\o=\i\ddbar f$ with the normalisations
$$\int_Me^{tf+c_t}\o^n=\int_Me^{f-(t-1)\vpt}\o^n=V.$$
Note that the first path is the
one used in Yau's continuity method proof \cite{Ya}. It connects any
point $\o$ in $\Hc$ to $\ricm1\o$ in $\hc2$. The second path, introduced by
Aubin in \cite{A2}, is used to connect any point in $\hc2$ to a \KE
metric.

The existence of the first path is equivalent to the Calabi-Yau theorem. The second path may not 
exist in the presence of nontrivial holomorphic vector fields but Bando and Mabuchi 
show that arbitrarily close to $\o$ in the $\CinfM$-topology
there exist metrics for which such a path exists. Since the K-energy
is continuous this will be sufficient for the argument (Cf. \cite{BM},\cite[Section 3]{SW}).

Now, for $t\in[0,1]$ one has 
\begin{equation}
 \Ric\ovpt=(1-t)\Ric\o+t\o,
\end{equation}
and
\begin{equation}
\D_t\dot\vpt=f+\dot c_t,
 \end{equation}
hence
\begin{eqnarray}
\ddt E_0(\oKE,\ovpt)
 &=&-\frac{1}{V}\intM\dot\vpt(\Ric\ovpt-\ovpt)\w n\ovpt^{n-1} \nonumber \\
 &=&-\frac{1}{V}\intM\dot\vpt((1-t)\i\ddbar f-\i\ddbar\ovpt)\w n\ovpt^{n-1} \nonumber\\
 &=&-(1-t)\frac{1}{V}\int_M\dot\vpt(\Dvpt\dot\vpt)^2\ovptn 
-\ddt(I-J)(\o,\ovpt) \nonumber\\
\label{FirstPatheq}
\end{eqnarray}
with $t\in[0,1]$, from which we conclude
\begin{equation}
E_0(\oKE,\o_{\vp_1})\leq E_0(\oKE,\o).\label{FirstEnergyDecreaseeq}
\end{equation}
Next, for $t\in[1,2]$
$$
\Ric\ovpt=(2-t)\o+(t-1)\ovpt,
$$
and
$$
\D_t\dot\vpt=-\vpt+t\dot\vpt,
$$
hence
\begin{eqnarray}
\ddt E_0(\oKE,\ovpt)
& =&-\frac{1}{V}\intM\dot\vpt(-(2-t)\i\ddbar\ovpt)\w n\ovpt^{n-1}\nonumber\\
& =&-(2-t)\ddt (I-J)(\oKE,\ovpt)\le0 \nonumber\\
& =&-(2-t)\frac{1}{V}\intM((\D_t\dot\vpt)^2+t|\del\dot\vpt|^2_t)\ovpt^{n}\le0,
\label{SecondPatheq}
\end{eqnarray}
where $t\in[1,2]$.
The theorem now follows in this case. Song and Weinkove extended this argument to $E_1$ using two detailed computations.
The first shows that while $E_k$ may not necessarily be monotone
(when the path exists), one still has $E_k(\oKE,\ovpone)\ge E_k(\oKE,\ovptwo)=0$.
In other words,
\begin{theorem-others} \label{SWThm}
\cite[Theorem 1.1]{SW} Let $(M,\oKE)$ be a Fano \KE manifold. 
Then for any $\o\in\hcbar2$ and for each $k=0,\ldots,n$
one has $$E_k(\oKE,\o)\ge0,$$ with equality if and only if $\o$ is
\KE and $h^\star\oKE=\o$ with $h$ a biholomorphism homotopic to the identity. 
\end{theorem-others}

\noindent The second calculation shows that when $k=1$, one has $$E_1(\ovpone,\o)\ge E_1(\o_0,\o)=0.$$
Explicitly, their computation shows that
\begin{eqnarray}
E_k(\ovpone,\o)=&
\frac{1}{V}\int_M
\i\del\vpone\w\dbar\vpone\w\sum_{i=0}^{n-1}a_i\o^{i}\w\ovpone^{n-1-i}\nonumber\\
&+(k+1)\frac{1}{V}\int_{M\times[0,1]}(1-t)(\Dvpt\dot\vpt)^2\ovptn\w dt \nonumber\\
&-\frac{1}{V}\int_M\sum_{i=1}^k \binom{i+1}{k+1}f(\i\ddbar
f)^i\w\o^{n-i},&\label{EkEnergyDecreaseeq}
\end{eqnarray} 
with
$a_i=\frac{(n-k)(i+1)}{n+1}$ if  $0\leq i\leq k-1$
and $a_i=\frac{(k+1)(n-i)}{n+1}$ if $k\leq i\leq n$. Since the last term is
positive on $\Hc$ for $k=1$ they conclude their proof. \bigskip

\subsection{A system of \MA equations}
Now, let $\o=\o_0$ denote an initial K\"ahler metric for our iterations. We present the iterative procedure defined by (\ref{iter_einstein}) in terms of \MA equations. Let $\vp_1$ be a \K potential
with 
$$
\Ric(\o_0+\i\ddbar\vp_1)=\o_0.
$$
Set the {\it Ricci deviation} $h:=h_{\o_0}$ of $\o_0$ as
$\i\ddbar h=\Ric\o_0-\o_0$. The function $h$ thus given is for the moment determined 
only up to a an additive constant.  The equation then becomes 
\begin{eqnarray}
-\i\ddbar\log\det g_{\vp_1}= & \o_0=-\i\ddbar\log\det g-\i\ddbar h
\end{eqnarray}
or
$$
\i\ddbar\log\frac{\o^n_{\vp_1}}{\o^n}=\i\ddbar h
$$
that is
$$
\o_{\vp_1}^n=e^h \o^n
$$
together with the volume normalisation 
$$
\V\int_M e^{h}\o^n=1.
$$
This determines $\vp_1$ only up to a constant, which will be fixed in the the next step.
Put $\o_1=\o_{\vp_1}$. In the second step we solve 
$$
\Ric(\o_1+\i\ddbar\vp_2)=\o_1
$$
and 
$\o_1-\Ric\o_1=\o_1-\o_0=\i\ddbar\vp_1$. 
The \MA equation is now 
$$
\o_{\vp_1+\vp_2}^n=e^{-\vp_1}\o_{\vp_1}^n= e^{h-\vp_1}\o^n,
$$
with $\vp_1$ determined uniquely by
$$
\V\intM e^{h-\vp_1} \o^n=1.
$$
Iterating this procedure we have 
$\ric{-l}\o=\o_{\Sigma^l_{j=1}\vp_j}$ for each $l\in\NN$
where
\begin{equation}
(\o+ \i\ddbar \Sigma^l_{j=1}\vp_j)^n =e^{h-\Sigma^{l-1}_{j=1}\vp_j}\o^n,\label{ConeposEq}
\end{equation}
and each of the $\vp_j$ is uniquely determined by
\begin{equation}
\V\intM e^{h-\Sigma^{l-1}_{j=1}\vp_j}\o^n=1.\label{ConeposnormalizationEq}
\end{equation}
From now on we set 
$$\Phi_l=\Sigma^l_{j=1}\vp_j$$ and 
$$\ol=\o_{\Phi_l}=\o_0+\i\ddbar \Phi_l.$$

\subsection{Monotonicity of the energy functionals}
The following proposition describes the monotonicity of the K-energy and $\Fmu$ along the iteration.
Note that the equivalent for the K\"ahler-Ricci flow is well-known.
\begin{proposition}
Let $\mu=1$. Then 
\begin{eqnarray}
E_0(\o,\ol) &=& -(I-J)(\o,\ol)-\V\intM
\Phi_{l-1}\oln
+\V\intM h\o^n
\leq 0, \nonumber\\
\Fmu(\o,\ol)&=& F^0(\o,\Phil) \le0,\nonumber\\
E_1(\o,\ol) &\leq & 0,\nonumber
\end{eqnarray}
with equality if and only if $\o$ is \KEno.
Furthermore all the $E_k$ with $k=2,\ldots,n$ decrease along the iteration starting from the second iteration.
\end{proposition}

\begin{proof} To show the first inequality we note
that
\begin{eqnarray}
E_0(\o_{k-1},\o_k)
&= & 
\V\intM-\vp_{k-1}\o_k^n-(I-J)(\o_{k-1},\o_k)
-\V\intM\vp_{k-1}(\o_{k-1}^n-\o^n_k) \nonumber\\
&= & 
-(I-J)(\o_{k-1},\o_k)-\V\intM\vp_{k-1}\o_{k-1}^n.
\end{eqnarray}
The first term is nonpositive with equality if and only if $\o_k=\o_{k-1}=\Ric\o_k$, while the second term
is nonpositive since
\begin{equation}\label{neg}
1=\V\intM \o_k^n=\V\intM e^{-\vp_{k-1}}\o_{k-1}^n\ge\V\intM(1-\vp_{k-1})\o_{k-1}^n.
\end{equation}
Since 
$$
E_0(\o,\ol) = \sum_{k=1}^lE_0(\o_{k-1},\o_k),
$$
the conclusion follows. \\
The second inequality follows similarly, indeed
$$
F_1(\o_{k-1},\o_k)=-(I-J)(\o_{k-1},\o_k)-\V\intM\vp_k\o_k^n,
$$
and with (\ref{neg}). The third inequality follows from
$$
E_1(\o_{k-1},\o_k)
=2F_1(\o_{k-1},\o_k)+\V\intm \vp_k(\o_k^n+\o_k^{n-1}\w\o_{k-1})
$$
from the formula relating $F_1$ and $E_k$. Since both summands are negative, we
conclude (note that $\o_{k-1}-\o_k=-\i\ddbar \vp_k$).
Finally, the decrease of the functionals are proved using (\ref{FirstEnergyDecreaseeq})
and (\ref{EkEnergyDecreaseeq}).
\end{proof}

\begin{remark}
For the rest of the energy functionals, 
we get by direct computations that with respect to a \KE metric $\o$,
\begin{eqnarray*}
E_k(\o,\ovp)&=&(k+1)F_1(\o,\ovp)-\V\intM\fovp\big(\ovpn+\ldots+\ovp^{n-k}\w(\Ric\ovp)^k\big)\\
&=&(k+1)F_1(\o,\ovp)-B_k(\o,\ovp).
\end{eqnarray*}
where $\Ric\ovp=\ovp+\i\ddbar\fovp$. Note that when the Ricci curvature is positive we have
\begin{eqnarray}
B_k(\o,\ovp)\hspace{-.25cm}&=& \hspace{-.25cm}
\V\intM\fovp(\ovp+\i\ddbar\fovp)\w(\ovp^{n-1}+\ldots+\ovp^{n-k}\w(\Ric\ovp)^{k-1})
\cr
&= &
-\V\intM\i\del\fovp\w\dbar\fovp\w(\ovp^{n-1}+\ldots+\ovp^{n-k}\w(\Ric\ovp)^{k-1})
\cr
&&\q +
\V\intM\fovp\ovp\w(\ovp^{n-1}+\ldots+\ovp^{n-k}\w(\Ric\ovp)^{k-1})
\cr
&\leq & B_{0}(\o,\ovp) \le0.
\end{eqnarray}
\end{remark}

\section{The dynamics of the Ricci operator}\label{SectionProofs}

We are now ready to answer the question raised by Nadel \cite{Na}. Note that some of the results presented in that section were discussed in details by the author and Y.\ Rubinstein  \cite{Ru} during the period\footnote{September 2003 to December 2006.} of their collaboration. 

\begin{theorem}\label{FirstNadelThm}
Let $(M,\o)$ be a \K manifold with positive first Chern class and assume that $\ric k(\o)=\o$
for some $k\in\NN$. Then $\o$ is \KEno.
\end{theorem}

\begin{proof}
Note that the nonexistence of fixed points of negative order implies that
of positive order, and reversely. Therefore assume that for some $\o\in\Hc$ and some
$l\in\NN$ one has $\ricm l(\o)=\o$. By the cocycle condition we therefore have
\begin{equation}
0=E_0(\o,\ricm l\o)=\sum_{i=0}^{l-1} E_0(\ricm{i}\o,\ricm {i-1}\o).\label{Ezerosumeq}
\end{equation}
On the other hand, from the first part of (\ref{MApathseq})
\begin{eqnarray*}
E_0(\o,\ricm 1\o)=
-\frac{1}{V}\int_{M\times[0,1]}(1-t)(\Dvpt\dot\vpt)^2\ovptn\w dt -(I-J)(\o,\ricm 1\o) 
\end{eqnarray*}
Thus $E_0(\o,\ricm 1\o)\leq 0$, with equality if and only if $\ricm 1\o=\o$. Therefore each of the terms in 
(\ref{Ezerosumeq}) must vanish identically and we conclude that $(M,\o)$ is \KEno. 
\end{proof}

\begin{proposition}\label{SecondNadelThm}
Let $M$ be as above and assume that $\o\in\hc k$ for all $k\in\NN$
and that $\o$ is not \KEno.
Then $\lim_{k\ra\infty}\ric k\o$ does not exist in $\hc0$.
\end{proposition}

\begin{proof}
If $\o_\infty=\lim_{l\ra\infty}\Ric^{(l)}\o$ exists and is 
smooth it satisfies $\Ric\o_\infty=\o_\infty$. But $E_0(\o,\o_\infty)>0$ contradicting
(\ref{BMLowerboundThm}).
\end{proof}

Let $G$ denote the Green function for $\D=\D_{\dbar}$ with respect to
$(M,\o)$ with $\intM G(x,y)\on(y)=0$ and $A(\o)=-\inf G$ such that
$$
f(x)-\VintM f\on = -\VintM G(x,y)\D f(y)\on(y),\q\all\, f\in\CinfM.
$$
Then, one has the following estimate due to Bando and Mabuchi.
\begin{theorem-others}\label{BMdiameterThm}
\cite{BM}  One has
$$
A(\o)\leq \frac {c_n}2\diam(M,\o)^2.
$$
If ${\rm Ric}(\o)\ge\eps\o$ for some $\eps>0$ then $\diam(M,\o)^2\le\frac{\pi^2(2n-1)}{\eps}$ by Myers' theorem. 
\end{theorem-others}
As an immediate corollary we have
\begin{lemma}\label{lemma1}
Let $M$ be a Fano manifold. Assume that the K-energy is proper. 
Let $(\o_l)_{l\in\NN}$ be  a sequence of \K forms on which the K-energy
is bounded from above and such that there 
exists  $l_0(\o_0)\in \NN$ and  $\eps>0$ 
with ${\rm Ric}\ol\ge\eps \ol,\all\,l\geq l_0$. Then
there exists a constant $C_1$ depending only on $(M,\o_0)$ such that
$$
\n{\Phi_l}\LinfMo\leq C_1,\q \all\,l\in\NN
$$
where $\o_l=\o_0+\i\ddbar \Phi_l$ and $\frac{1}{V}\int_M e^{h-\Phi_l}\frac{\o_0^n}{n!}=1$.
\end{lemma}

\begin{proof} 
Let $G_l$ be the Green function for $\D_l=\D_{\dbar,\ol}$ 
(i.e., the Laplacian with respect to $(M,\ol)$) 
satisfying $\intM G_l(x,y)\oln(y)=0$. Set $A_l=-\inf_{M\times M} G_l$.

Since $-n<\D_1\Phil$ and $n>\D_l\Phil$ the Green formula gives
\begin{eqnarray}
\Phil(x)-\VintM \Phil\o_0^n = & -\VintM G_1(x,y)\D\Phil(y)\o_0^n(y)\leq nA_1,\\
\Phil(x)-\VintM \Phil\oln= & -\VintM G_l(x,y)\D\Phil(y)\oln(y)\geq -nA_l.
\end{eqnarray}
Hence
\begin{equation}
\osc\Phil\leq n(A_1+A_l)+I(\o,\ol).\label{OscEq}
\end{equation}

Since $E_0$ is proper on $\HcG$ in the sense of Tian, if $E_0(\o_0,\cdot)$ is uniformly bounded from above on a subset of $\HcG$ so is $I(\o,\cdot)$. 
We conclude that $I(\o_0,\ol)$ is uniformly bounded independently of $l$.
Finally, the Proposition follows from Theorem \ref{BMdiameterThm} which provides a uniform 
bound for $A_l$. 
\end{proof}

As a consequence of the properness of the K-energy for Fano Einstein manifolds \cite{T2,PSSW} we obtain the following,
 
\begin{corollary}
Let $M$ be a Fano Einstein manifold with no nontrivial holomorphic vector field. Consider the sequence of K\"ahler metrics $(\o_l)_{l\in \NN}$ defined by the system of \MA equations (\ref{ConeposEq}) and assume that for $l$ sufficiently large there exists a constant $\epsilon>0$ with 
$$\Ric(\o_l) \geq \epsilon \o_0.$$
Then $\o_l$ converges to the \KE metric when $l$ tends to infinity.   
\end{corollary}

Finally, inspired by \cite{Pa}, we derive the following,
\begin{proposition}
Let $M$ be a Fano manifold with $G$ a maximal compact subgroup of $\AutMJ$. Consider the sequence of $G$-invariant K\"ahler metrics $\o_l$ defined by the system of \MA equations (\ref{ConeposEq}). Assume that there exists a constant $1>\kappa>0$ such that for $l$ sufficiently large,
$$(2-\kappa)\o_0^n \geq \o_l^n \geq \kappa \o_0^n$$
where $\o_0$ is $G$-invariant K\"ahler metric. Then $M$ is \KE and $\o_l$ converges to a $G$-invariant \KE metric when $l$ tends to infinity.
\end{proposition}
\begin{proof}
Thanks to the proof of Lemma \ref{lemma1} and Theorem \ref{BMdiameterThm}, we are reduced
to prove an upper bound for $I(\o,\ol)$.
But if we denote $\Phi_l^+(x)=\sup\{0,\Phi_{l}(x)\}$ and $\Phi_l^-(x)=\inf\{0,\Phi_{l}(x)\}$, we obtain
\begin{eqnarray*}
I(\o,\ol) &\leq&\frac{1}{V}\int_{\{\o_l^n \geq \o_0^n\}} (-\Phi_l^-)(\o_l^n - \o_0^n) +
\frac{1}{V}\int_{\{\o_0^n \geq \o_l^n\}} (\Phi_l^+)(\o_0^n - \o_l^n) \\
&\leq& (1-\kappa)\frac{1}{V}\int_M (\Phi_l^+ - \Phi_l^- ) \o_0^n\\
&\leq& (1- \kappa)\osc\Phil
\end{eqnarray*}
Together with 
\begin{equation*}
\osc\Phil\leq n(A_1+A_l)+I(\o,\ol).
\end{equation*}
this gives us to the $\C^0$ bound for $\Phi_l$. Now it is a standard argument of \MA equations to get the convergence in the $\C^{\infty}$ topology.
\end{proof}

For a K\"ahler manifold $M$, we consider the family of \MAes \- (\ref{MApathseq}).
We introduce the {\it Aubin operators} $\Aubeps$ by setting 
$$\Aubeps(\o)=\o_{\vp_{1+\eps}}$$ for each $\eps\in[0,1]$ such that
$\vp_{1+\eps}$, solution of (\ref{MApathseq}), exists. Note that 
$\Aub_0(\o)=\ricm1\o$ and $\Aub_1(\o)=\oKE.$ Formally, one can think
$$\Aubeps=\left(\frac1{1-\eps}(\Ric-\eps\II)\right)^{-1}$$
and that we have defined the following sequence of \MA equations
\begin{equation}
(\o_0 + \i\ddbar \Phi_{j})^n=e^{(\eps-1)\Phi_{j-1}-\eps \Phi_j}(\o_0+\i\ddbar \Phi_{j-1})^n
\end{equation}
Let $G\sseq\AutMJ$ be as before. Then from uniqueness of solutions of the family of the \MA equations (\ref{MApathseq}) we conclude that
$\Aubeps$ maps $\HcG$ into itself.

We recall the definition of $\a$-invariant introduced by Tian.
\begin{equation}
\a_{M}=\sup\{\a\ge0\,:\, \sup_{\vp\in\Hc}\int_M e^{-\a\vp}\o^n<\infty\}.\label{Alphaeq}
\end{equation}
In \cite{T0} it is proved that $\a_{M}$ is a positive
holomorphic invariant of Fano manifolds. Regarding the existence of the Aubin operators we state the following
\begin{proposition}\-\\
(i) \cite{T0} Assume that $M$ is Fano. Then the operators 
$\Aubeps$ exist for all
 $\eps\in[0,\min {\{1,{\frac{n+1}{n}}\alpha_{M}\}})$. 
\\(ii) \cite[Theorem 5.7]{BM} Assume in addition
that the K-energy is bounded from below. Then the operators $\Aubeps$ exist for
any $\eps\in[0,1)$. 
\end{proposition}

We now recover by a conceptually simpler method a theorem of Tian \cite{T2}.

\begin{corollary}
Let $M$ be a Fano manifold. Let $G$ be a maximal compact subgroup of $\AutMJ$ and assume that the K-energy is proper on $\HcG$ and let $\eps\in(0,1)$. Then there exist $G$-invariant \KE metrics. All such metrics 
are the limit points of the iterates of $\Aubeps$ on $\HcG$ in the $\Cinf(M)$-topology.
\end{corollary}

\begin{proof}
Since the K-energy is proper on $\HcG$ and in particular bounded from below, $\Aubeps$ are defined
for each $\eps\in(0,1)$. As before we obtain (\ref{OscEq}) with $\ol=\Aubepsl(\o)$ with each
$\ol$ in $\HcG$. 
Since 
by (\ref{FirstEnergyDecreaseeq}) and (\ref{SecondPatheq})
the K-energy decreases along iterates we still have a uniform bound on
$I$ along the orbits. By Theorem \ref{BMdiameterThm} we also have a uniform bound (depending on $\eps$) on 
$A(\Aubepsl(\o))$. 
For each $\be\in(0,1)$ the $\Ctwobe(M)$-estimates now follow from \cite[Chapter 7]{A4} and the higher order derivatives follow from these, by bootstrapping. We may therefore extract
a converging subsequence $\Aubepslj(\o)$ in the $\CtwobeMo$-topology whose limit $\o_\infty$ lies in $\CtwoM$. Moreover, since the 
K-energy is bounded from below we have 
$$
\lim_{j\ra\infty}E_0(\Aubepsj{l_j}(\o),\Aubepsj{l_{j+1}}(\o))=
\lim_{j\ra\infty}\sum_{k=l_j}^{l_{j+1}-1}E_0(\Aubepsj{k}(\o),\Aubepsj{k+1}(\o))=0
.$$ 
As each of the summands is nonpositive one has
$E_0(\o_\infty,\Aubeps(\o_\infty))=0$ and it follows that $\o_\infty\in\HcG$ is smooth and \KE  
by
(\ref{FirstEnergyDecreaseeq}) and (\ref{SecondPatheq}). Since this is true for each converging subsequence we conclude that, in fact, the sequence of iterates itself converges.
\end{proof}

In fact we also get,

\begin{theorem}\label{AubepsThm}
Let $M$ be a \KE manifold with positive first Chern class. Let
$G$ be a maximal compact subgroup of $\AutMJ$.

Then for any $\eps\in(0,1)$, $\o\in\HcG$ there exists a  
biholomorphism $h$ of $(M,J)$ such that one has 
$\lim_{l\ra\infty}\Aubepsl\o=h^\star\oKE$ in the $\CinfM$-topology for some $G$-invariant
\KE form $\oKE$. 
\end{theorem}

\section{The case of toric Fano manifolds}\label{toric}

On a toric K\"ahler manifold $M$ with positive first Chern class, a result of X-J.\ Wang and X.\ Zhu \cite{WZ} asserts that there does always exist a K\"ahler-Ricci soliton, which is unique up to holomorphic automorphims. A K\"ahler-Ricci soliton is a pair $(X,\omega_g)$ where $X$ is a holomorphic vector field on $M$ and $\omega=\omega_g$ a K\"ahler form that satisfies

\begin{equation}\label{Soliton}
\Ric(\omega)=\omega+ L_X(\omega)
\end{equation}
where $L_X$ is the Lie derivative along $X$. In fact, the K\"ahler-Ricci soliton in the toric case is a genuine K\"ahler-Einstein metric if and only if the Futaki invariant vanishes which implies that the holomorphic automorphism group of the manifold is reductive.
\par Given a nontrivial holomorphic vector field $X$ on $M$, it is well known by Hodge theory that there exists a unique smooth complex-valued function $\theta_X$ such that
$$i_X(\omega)=\i\bar\partial \theta_X$$
with $\int_M e^{\theta_X} \omega^n = V$. Let $h\in \CinfM$ be the Ricci deviation of $\omega$ determined by $\frac{1}{V}\int_M e^h \omega^n=1$. From \cite[Lemma 2.1]{WZ}, there exists a unique holomorphic vector field $X$ such that $$F_X(v):=\int_M v(h-\theta_X)e^{\theta_X}\omega^n=0$$
for all holomorphic vector field $v$. Note that $F_X$ is actually independent of the choice of $\omega$. 
For such a vector field $X$, to find a K\"ahler-Ricci soliton $(X,\omega+\i\ddbar \varphi)$ satisfying (\ref{Soliton}) is equivalent to solve the \MA equation
$$(\omega+\i\ddbar \varphi)^n = e^{h-\theta_X - X(\varphi) -\varphi}\omega^n$$
In \cite{Zh}, an existence result (proved by a continuity method argument) is given for the \MA equations of the form $$\text{det}(g_{i\overline j}+\varphi_{i\overline j})=e^{f-t\theta_X-tX(\varphi))}.$$ where $f$ is smooth and $t$ varies from 0 to 1 (note that $t=0$ corresponds exactly to Calabi-Yau theorem). This can be seen as consequence of the work of S.\ Kolodziej \cite{Ko} for solving \MA equations, and the following result \cite[Corollary 5.3]{Zh}.
 \begin{lemma}
  Let $(M,\omega)$ a K\"ahler manifold with a non-trivial holomorphic vector field $X$. Suppose that $\varphi$ is smooth, $\omega + \sqrt{-1}\ddbar \varphi >0$ and $X(\varphi)$ is a real-valued function. Then there is a uniform constant $C$ independent of $\varphi$ such 
  that $|X(\varphi)|<C.$ 
 \end{lemma}
  Hence, we can define naturally the sequence of \MA equations
 $$(\omega_{j-1}+\i\ddbar \varphi_j)^n = e^{h-\theta_X - X(\varphi_{j}) -\varphi_{j-1}}\omega_{j-1}^n$$
which turns out to be equivalent to solve for all $j\geq 1$
\begin{equation}\label{iter_soliton}
\Ric(\omega_j)=\omega_{j-1}+L_X(\omega_{j})
\end{equation}
for $\omega_j=\omega_{j-1}+\i\ddbar \varphi_j=\omega_0+ \i\ddbar \Phi_j$. \bigskip
\par The goal of this section is to prove that under this setting, one can always get a $\C^0$ estimate for the potentials $\Phi_j$, which means that the iteration procedure defined by (\ref{iter_soliton}) converges to the unique K\"ahler-Ricci soliton on the manifold.
We follow the techniques developed in \cite{WZ} using real \MA equations to derive the uniform estimate for $\Phi_j$.
\par We shall need the following lemma that we will apply with the convex polyhedron $\Omega^*$ associated to our toric manifold.
\begin{lemma} \label{lemma-center}
Let $\Omega$ a bounded convex domain in $\mathbb{R}^n$. Then there is a unique ellipsoid $\mathcal{E}$ called the minimum ellipsoid of $\Omega$ which attains minimum volume among all ellipsoids containing $\Omega$, such that
$$\frac{1}{n} \mathcal{E} \subset \Omega \subset \mathcal{E} $$
where $r \mathcal{E}$ denotes the dilatation of factor $r$ of $\mathcal{E}$ with concentrated factor.
\end{lemma}
Note that the vector field $X$ can be expressed in the form 
$$X=\sum_i {c}_i X_i$$ 
with $X_i=\frac{\partial}{\partial \gamma_i}$ for $(\gamma_i=x_i+\i\theta_i)_{i=1,..n}$  affine logarithm coordinates on the maximal torus $T$ of $\AutMJ$. In \cite[Section 2]{WZ} it is shown that these constants ${c}_i$ satisfy the relations
\begin{equation}\label{relation_coeff}
\int_{\Omega^*} y_i e^{\sum_{l=1}^n c_l y_l } dy =0 
\end{equation}
for all $i=1,..,n.$. Furthermore, there exists a K\"ahler form $\omega_0$, invariant under the maximal compact abelian subgroup $K_0$ of $T$, and determined by a convex function $w_0$ on $\RR^n$, i.e on $T$
\begin{equation} \label{start_iter_soliton}
\omega_0 = \i \ddbar w_0.
\end{equation}
In all the following, we shall denote $\det(w_0):=\det\left ( (w_0)_{ij} \right)$. We obtain
$$\omega_0^n =\frac{1}{\pi^n} \det(w_0) dx_1 \wedge ...\wedge dx_n \wedge d\theta_1 \wedge ... \wedge d\theta_n.$$
Note that the function $w_0$ is determined uniquely by the convex polyhedron $\Omega^*$. Indeed, if $\{p^{(i)}\}_{i=1,..,m}$ denote the vertices of $\Omega^*$, $w_0$ is the convex function defined by $w_0=\log \left( \sum_{i=1}^m e^{\langle p_i ,x \rangle } \right) $. 
After normalisation and without loss of generality, we can assume $\det(w_0)=e^{-h-w_0}$. We initialize our iteration procedure by chosing this way the K\"ahler form $\o_0$. Let's call $w_j=w_0+\Phi_j$ with $\Phi_j$ a $K_0$-invariant potential. Now, there exists a constant $c_X$ such that
$$\theta_X + X(\Phi_j)= \sum_{l=1}^n c_l \frac{\partial w_j}{\partial x_l} + c_X.$$
This means that to solve (\ref{iter_soliton}) is equivalent to solve the real \MA equation on $\RR^n$,
\begin{equation}\label{real_iter_soliton}
\det(w_0+\Phi_j) =\det(w_j)= e^{-c_X - w_{j-1} - \sum c_l \frac{\partial w_{j}}{\partial x_l}}
\end{equation}

\begin{lemma}
Define
$$m_j=\inf_{\RR^n} w_j.$$ 
Then there exists a constant $C$ independent of $j$ such that 
 $$m_j \leq C.$$
\end{lemma}
\begin{proof}
For an integer $k$, we choose $x_{j,k}$ a point in $\RR^n$ such that
$w_{j}(x_{j,k})=m_{j}+k$. We introduce for $j \geq 0$, the set
\begin{eqnarray*}
A_{k,j}=\lbrace x \in \RR^n &:& m_{j}+k \leq w_{j}(x) \leq m_{j}+k+1,  \\
            && m_{j+1}+r_{j+1,k} \leq w_{j+1}(x) \leq m_{j+1}+r_{j+1,k}+1 \rbrace 
\end{eqnarray*}
where $r_{j+1,k}=[w_{j+1}(x_{j,k})-m_{j+1}]+k$. 
For any integer $k \geq 0$, $\cup_{0\leq i\leq k}^{k} A_{i,j}$ is convex. On the other hand, we know that $D\omega_j(\RR^n)=\Omega^*$ and thus $A_{k,j}$ is bounded.
Let's define $d_{\max}=\sup_{y \in \Omega^*} \{ c_l y_l \}$. We get on ${A}_{0,j}$,
$$\det(w_{j+1}) \geq e^{-c_X-d_{\max}-1} e^{-m_{j}}.$$
From Lemma \ref{lemma-center}, there exists a linear transformation $y=Tr(x)$ which leaves the center of the minimum ellipsoid of ${A}_{0,j}$ invariant such that
$B_{R/n} \subset Tr({A}_{0,j}) \subset B_R$ for a certain well-chosen constant $R$. Let's try to obtain some information about this radius. 
Define $$v_j=\frac{1}{2} \left(e^{-c_X-d_{\max}-1}\right)^{\frac{1}{n}} e^{-\frac{m_{j}}{n}} \left( |y-y_j|^2- \left(\frac{R}{n}\right)^2 \right) + m_{j+1}+r_{j+1,0}+1 $$
where $y_j$ the center of the minimum ellipsoid of ${A}_{0,j}$. Then, in $Tr({A}_{0,j})$,
$$\det(v_{j}) = e^{-c_X-d_{\max}-1}e^{-m_{j}}$$
and of course, $v_j \geq w_{j+1}$ on $\partial Tr({A}_{0,j})$. The comparison principle for \MA operator gives that $v_j \geq w_{j+1}$ on $Tr({A}_{0,j})$. Hence,
\begin{eqnarray*}
m_{j+1}+r_{j+1,0} &\leq& w_{j+1}(y_j)\\ 
                   &\leq& v_j(y_j)\\
                   & =& -\frac{1}{2}  \left(e^{-c_X-d_{\max}-1}\right)^{\frac{1}{n}}e^{-\frac{m_{j}}{n}}\left(\frac{R}{n}\right)^2 + m_{j+1}+r_{j+1,0} +1
\end{eqnarray*}
This gives us the following inequality with $C_0=\sqrt{2}n\left(e^{-c_X-d_{\max}-1}\right)^{-\frac{1}{2n}}$,
$$R\leq  C_0 e^{\frac{m_{j}} {2n}}$$
On the other hand, since $Tr(A_{k,j}) \subset B_{2(k+1)R}$, we obtain
\begin{eqnarray*}
\int_{\RR^n} e^{-w_{j}} &=& \sum_{k} \int_{Tr(A_{k,j})} e^{-w_{j}}\\
 &\leq  & \sum_{k} e^{-m_{j}-k} Vol(Tr(A_{k,j}))\\
 & \leq & Vol(S^{n-1}) \sum_{k} e^{-m_{j}-k} 2(k+1)R^n\\
 & \leq & C_{n} e^{-m_{j}/2}
\end{eqnarray*}
But, on another hand, we also have thanks to (\ref{real_iter_soliton}) and by transformation $y=Dw_{j+1}(x)$,
\begin{eqnarray}\label{boundedness_toric}
\int_{\RR^n} e^{-w_{j}} dx= e^{c_X}\int_{\Omega^*} e^{\sum_l c_l y_l} dy = C_0.
\end{eqnarray}
This gives the expected inequality.
\end{proof}
 
\begin{lemma}\label{upper_bound_x_j}
Let $x_{j,0}$ be the minimum point of $w_{j}$. Then there exists $C>0$ such that $|x_{j,0}|<C$ for all $j$.
\end{lemma}
\begin{proof}
The fact that $|Dw_{j}|\leq \sup_{x \in \Omega^{*}}|x|$, and the previous lemma
gives the existence of $R'$ independent of $j$ such that $\inf_{\partial B_{R'}(x_{j,0})} w_{j} \leq m_{j} +1$. By convexity, 
\begin{equation} \label{ineq0}
|Dw_{j}(x)|\geq 1/R'
\end{equation}
in $\RR^n \backslash B_{R'}(x_{j,0})$.
Hence with (\ref{boundedness_toric}), for any $\epsilon>0$, there exists $R_\epsilon$ sufficiently large such that
$$\int_{\RR^n \backslash B_{R_\epsilon}(x_{j,0})} e^{-w_{j}} dx \leq C_1 \int_{\RR^n \backslash B_{R_\epsilon}(x_{j,0})} e^{|x-x_{j,0}|/R'} < \epsilon$$
where $R_\epsilon$ is independent of $j$.\par
Now, assume that $|x_{j,0}|$ is not bounded. For any $\epsilon>0$, there exists a large constant $C>0$ such that if $|x_{j,0}|>C$ then, $$\frac{\partial w_{j+1}}{\partial \xi} >\frac{1}{2} \inf_{x \in \partial\Omega^*}{|x|}$$ for $x\in B_{R_\epsilon}(x_{j,0})$ and where we have set $\xi=\frac{x_{j,0}}{|x_{j,0}|}$. This can be viewed by considering the restriction of $w_{j+1}$ on the ray $\overrightarrow{Ox_{j,0}}$ and by its convexity together with the fact that $Dw_{j+1}(\RR^n)=\Omega^*\ni \{0\}$. Now with (\ref{boundedness_toric}),
\begin{eqnarray}\label{ineq1}
\int_{B_{R_\epsilon}(x_{j,0})} \frac{\partial w_{j+1}}{\partial \xi} e^{-w_{j}} \geq \frac{C_0}{4}\inf_{x \in \partial\Omega^*}{|x|}.
\end{eqnarray}
Furthermore, we have
\begin{equation} \label{ineq2}
\left\vert \int_{\RR^n \backslash B_{R_\epsilon}(x_{j,0})} \frac{\partial w_{j+1}}{\partial \xi} e^{-w_{j}} \right\vert \leq  \left(\sup_{x \in \Omega^{*}}|x|\right) \left\vert{\int_{\RR^n \backslash B_{R_\epsilon}(x_{j,0})} e^{-w_{j}}}\right\vert \leq   \epsilon \sup_{x \in \Omega^{*}}|x|
\end{equation}
With (\ref{ineq1}) and (\ref{ineq2}), one obtains for $\epsilon$ sufficiently small that
\begin{equation}
\int_{\RR^n}   \frac{\partial w_{j+1}}{\partial \xi} e^{-w_{j}} dx > 0.\label{contradiction}
\end{equation}
Now on another hand by (\ref{real_iter_soliton}) and (\ref{relation_coeff}), 
\begin{eqnarray*}
0&=&\int_{\Omega^*} y_i e^{\sum_{l=1}^n c_l y_l } dy \\
&=&\int_{\RR^n} \frac{\partial w_{j+1}}{\partial x_i} e^{\sum_l c_l \frac{\partial w_{j+1}}{\partial x_l} }\det(w_{j+1}) dx \\
 &=& e^{-c_X} \int_{\RR^n} \frac{\partial w_{j+1}}{\partial x_i} e^{-w_{j}} dx
\end{eqnarray*}
 which leads us to $\int_{\RR^n} \frac{\partial w_{j+1}}{\partial \xi} e^{-w_{j}}=0$ and hence a contradiction with (\ref{contradiction}). Thus, $|x_{j,0}|$ has to be bounded.
\end{proof}

\begin{lemma}
 One has the upper bound
$$\sup_{M} \Phi_j < C'$$
for a constant $C'$ independent of $j$.
\end{lemma}
\begin{proof}
As explained in the proof of \cite[Lemma 3.4]{WZ}, by convexity, it is sufficient to give an upper bound of $\Phi_j(0)$. Since $|x_{j,0}|<C$ by Lemma \ref{upper_bound_x_j} and $|Dw_j|\leq \sup_{x \in \Omega^*} |x|$, we obtain $w_j(0)<C'$. With (\ref{boundedness_toric}), we get $|w_j(0)|$ bounded and hence $|\Phi_j(0)|$ bounded.
\end{proof}

\begin{lemma}
One has the lower bound
$$\inf_{M} \Phi_j > C''$$
for a constant $C''$ independent of $j$.
\end{lemma}
\begin{proof}
Let $\{p^{(i)}\}_{i=1,..,m}$ denote the vertices of $\Omega^*$ and define $$\bar v(x)=\max_{i=1,..,m}\{ \langle x,p^{(i)} \rangle \}.$$  Then the graph of $\bar v$ is a convex cone with vertex at the origin and is an asymptotical cone of the graph of $w_0$, i.e one can check that the following inequality is satisfied $$|\bar v  - w_0|\leq c_0.$$
Let's denote $$z_j(r)=\sup_{|x|=r}( \bar v-w_j )(x)$$ and assume that the supremum is attained at $p_{r,j}$. Then with $\xi=\frac{p_{r,j}}{|p_{r,j}|}$, we get
\begin{equation}\label{ineq3}
 z_j'(r) \leq \partial_{\xi} (\bar v - w_j)(p_{r,j}). 
\end{equation}
 
\begin{claim}
There exists $r_0$ independent of $j$ such that for  $r>r_0$, $z_j'(r)\leq A r^{-2}$
for a constant $A$ independent of $j$. 
\end{claim}
If the claim is proved then by integration, $\sup_{x\in \RR^n}(\bar v -w_j)$ is bounded independently of $j$ and the lemma is proved.
Let's now prove the claim. Let $F^{(k)}$ be the faces of the graph of $\bar v$ and suppose that \begin{equation} \label{inclusion_F}
F^{(1)} \subset \{x_1 > c_0'|\hat x|\} \cap \{x_{n+1}=0\}
\end{equation} 
which is possible by a change of variables ($c_0'$ is here a well-chosen positive constant and $\hat x=(x_2,...,x_n)$). As explained in \cite[Lemma 3.5]{WZ}, if the claim does not hold, there would exist a point $q_{r,j}$ in $\{x_1=\frac{1}{2} p_{r,j}\}$ such that
$$w_j(q_{r,j})=\inf \{w_j(x) : x\in F^{(1)} \cap \{x_1=\frac{1}{2} p_{r,j}\} \}.$$
It follows from (\ref{ineq3}) that for a certain constant $c_0''$,
\begin{equation} \label{ineq4}
w_j(q_{r,j})-w_j(p_{r,j}) \geq \frac{c_0''}{r}.
\end{equation} 
Define for a convex function $f$ the set $$N_{f}(x)= \{p \in \RR^n : f(y)\geq f(x) + p\cdot (y-x) 
\hspace{0.35cm} \forall y\in \RR^n\}$$
and for a set $S$, $N_f(S)=\cup_{x \in S} N_f(x)$. Then, for $\psi$ defined on $\RR^n$ such that its graph is a convex cone with $\psi(p_{r,j})=w_j(p_{r,j})$ and $\psi=w_j$ on 
$\{x \in \RR^n: w_j(x)=w_j(q_{r,j})\}$, one has
\begin{equation} \label{inclusion}
N_{\psi}(S_0)\subset N_{w_j}(S_0) 
\end{equation}
where we have defined $S_0=\{x \in \RR^n: w_j(x)<w_j(q_{r,j})\}$. From (\ref{ineq4}), one gets $${\rm dist}(0,S_0) \geq c_2r $$
and by (\ref{ineq0}) applied at step $j-1$, we obtain for $x\in S_0$ the existence of a constant $c_3>0$ such that
$$w_{j-1}(x) \geq c_3 |x| - C \geq c_2 r-C.$$
But now, with (\ref{real_iter_soliton}), there exists a positive constant $c_3'$ such that
$$Vol(N_{w_j}(S_0))\leq \int_{S_0} \det(w_j) \leq c_3Vol(S_0)e^{-c_3'r}$$
and on another hand, with (\ref{inclusion_F}) and (\ref{ineq4}), there exists $c_4$ such that 
$$Vol(N_{\psi}(S_0)) \geq \frac{c_4}{r^{n+1}}$$
which leads us to the expected contradiction because of the inclusion (\ref{inclusion}).
\end{proof}

Once the $\C^0$ estimate in hand for $\Phi_j$, one can derive $\C^\infty$ estimates by using the results of \cite{Ya,TZ}. Finally the results of this section gives the
\begin{theorem}
Let $M$ be a Fano toric manifold with K\"ahler-Ricci soliton $(X,\omega_{\infty})$. The iterations defined by   
$$Ric(\omega_j)-L_X(\omega_j)=\omega_{j-1} $$
for $j\geq 1$ and with $\omega_0$ given by (\ref{start_iter_soliton}), converge to
a K\"ahler-Ricci soliton $(X,\omega_\infty)$. If the Futaki invariant vanishes, then the iterations converge to the K\"ahler-Einstein metric on $M$.
\end{theorem}
\begin{remark}
In fact our theorem gives another proof of the existence of K\"ahler-Ricci solitons on Fano toric manifolds without using any flow or continuity method.
\end{remark}
\section{Applications and numerical results}

\subsection{The finite dimensional picture}\label{finite}

In \cite{D3}, Donaldson has introduced the notion of $\nu$-balanced metric for a fixed volume form $\nu$ and proved its existence \cite[Prop. 4]{D3} under some very general conditions \cite[p.10]{D3}. These metrics have the properties to solve the Calabi problem, i.e to converge towards the K\"ahler metric that has volume form $\nu$ in a given K\"ahler class.
\par Let us fix a volume form $\nu$ on a smooth projective manifold $M$ with a polarisation $L$, and
choose $r\in \NN$ sufficiently large such that $M$ is embedded by the holomorphic sections of $L^{r}$ in the projective space $\mathbb{P}H^0(M,L^r)$. We set $N_r=h^0(M,L^r)$ which is finite since $M$ is compact. One defines a $\nu$-balanced metric at rank $r$ as the fixed point of the map $T_\nu : Met(H^0(M,L^r)) \rightarrow Met(H^0(M,L^r))$,  
$$ T(G)_{i,j} = \frac{N_r}{Vol_L(M)}\int_M \frac{\langle S_i, S_j \rangle}{\sum_i |S_i|^2} d\nu$$
where $G$ is a hermitian metric of $H^0(M,L^r)$ and $(S_i)_{i=1,..,N_r}$ is an orthonormal basis of $H^0(M,L^r)$ with respect to $G$. Donaldson proved that the compositions of the map $T_\nu $ give a convergent sequence of metrics in $Met(H^0(M,L^r))$ (and thus on $Met(L^r)$ by the Fubini map $FS$, see \cite[p.4, Section 2.2]{D3}). The limit is called the $\nu$-balanced metric.
\begin{notation}
For a smooth hermitian metric $h\in Met(L)$ on the line bundle $L$, we denote $c_1(h)\in 2\pi [L]$ its curvature.
\end{notation}

\begin{theorem}
Under these settings, let's call $H_r\in Met(H^0(M,L^r))$ the sequence of $\nu$-balanced metrics of order $r$. Then $c_1(FS(H_r)^{1/r})$ converge to a K\"ahler form $\omega_{\infty}$ in $[c_1(L)]$ that satisfies
$$\omega_{\infty}^n = \nu.$$ 
\end{theorem}
\begin{proof}
To prove this theorem, we use the powerful Calabi-Yau theorem. Hence we know the existence of a K\"ahler form $\omega$ in $[c_1(L)]$ such that $\omega^n=\nu$. We use X.\ Wang's theorem \cite{W2} with the trivial bundle and $L$. There is a Hermitian-Einstein metric on these bundles and the metrics $H_r$ are \lq balanced\rq  \- with respect to $\omega$ in the sense studied by Wang. This is due to the obvious fact that the considered bundles are Gieseker stable. Thus, one obtains directly the convergence of the sequence of metrics
$FS(H_r)^{1/r}\in Met(L)$ to the metric $h_L$ with $c_1(h_L)=\omega$.
\end{proof}

\par We now assume in this section that $M$ is Fano and consider the polarisation $L=-K_M$. We fix an integer $r$ sufficiently large such that $M$ is embedded by the holomorphic sections of $L^{r}$. Let's consider a smooth hermitian metric $h_0$ on $L$ with $c_1(h_0)=\omega_0$ and let's call $f_{\omega_0}$ the Ricci deviation of $\omega_0$. Now for each $k$, one can define a $\Ric^{(-k)}(\omega_0)^n$-balanced metric at rank $r$ in the following way.

Consider $Hilb_{\omega_0}(h_0^r)$ the $L^2$-metric induced on the space of holomorphic sections $H^0(M,L^r)$ with respect to $\omega_0$, given by
$$Hilb_{\omega_0}\left(\langle,\rangle\right)(S_i,S_j)=\int_M \langle S_i,S_j \rangle e^{f_{\omega_0}}\omega_0^n.$$
Now, for a given hermitian metric $H_0$ on $H^0(M,L^r)$,  we define the metric $FS(H_0)$ on $L^r$ by
$$\sum_{i=1}^{N_r} |S_i|^2_{FS(H_0)}=\frac{N_r}{V} $$
where the $(S_i)_{i=1,..,N_r}$ form an $H_0$-orthonormal basis of $H^0(M,L^r)$.
Then, from \cite[Prop. 4]{D3}, we know that the dynamical system $FS \circ Hilb_{\omega_0} $ has an attractive fixed point $h_{\omega_0,r}$ at rank $r$ and that the convergence of this dynamical system is exponentially fast. We obtain this way a new form $$\omega_{1,r}=c_1(h_{\omega_0,r}).$$ For the second step, i.e in order to find the balanced metric $h_{\Ric^{-1}(\omega_0),r}$, we introduce the operator $Hilb_{\omega_{1,r}}$, by
$$Hilb_{\omega_{1,r}}\left(\langle,\rangle\right)(S_i,S_j)=\int_M \langle S_i,S_j \rangle e^{-\varphi_{0,r}}\omega_{1,r}^n$$
where $\varphi_{0,r}$ is the potential of the metric $h_{\omega_0,r}\in Met(L^r)$. Iterating this procedure leads us to define at each step $k$ a dynamical system $Hilb_{\omega_{k,r}} \circ FS$ which has an attractive fixed point $h_{\omega_{k,r}}\in Met(L^r)$ ($=e^{-\varphi_{k,r}}|.|_0$ locally). Note that for a generic $k$, we have
\begin{equation}\label{presque-donaldson}
Hilb_{\omega_{k,r}}\left(\langle,\rangle\right)(S_i,S_j)=\int_M \langle S_i,S_j \rangle e^{-\varphi_{k-1,r}}\omega_{k,r}^n
\end{equation}

\begin{corollary}
Under above assumptions, and for $r$ sufficiently large, the sequence $c_1(h_{\omega_{l,r}}^{1/r})$ converges when $r$ tends to infinity to the solution of the \MA equation (\ref{ConeposEq})  with exponential speed of convergence.
\end{corollary}
\begin{conjecture}
Under above assumptions, the sequence $c_1(h_{\omega_{l,r}}^{1/r})$ converges when $l$ tends to infinity to a K\"ahler metric $\omega_{r}$ in the class $c_1(M)$ with exponential speed of convergence. If $M$ is \KE, then $\o_r$ converges to a \KE metric.
\end{conjecture}

Let's describe now how our discussion can be useful for numerical approximations of \KE metrics on Fano manifolds. One has to notice at this stage that we can write  (\ref{presque-donaldson}) as
\begin{equation}\label{presque-donaldson2}
Hilb_{\omega_{k,r}}\left(\langle,\rangle\right)(S_i,S_j)=\int_M \langle S_i,S_j \rangle 
\left(\frac{V}{N_r} \sum_{i} \vert \widetilde{S_{i,k-1}} \vert^2 \right)^{-1/r} 
\end{equation}
where the $(\widetilde{S_{i,k-1}}) \in H^0(M,L^r)$ form an orthonormal basis of holomorphic sections with respect to the $L^2$ metric $H_{k-1,r}=Hilb_{\o_{k-1,r}} (h_{\omega_{k-1,r}})$ computed at the previous step, i.e the $\Ric^{-(k-1)}(\o_0)^n$-balanced metric. Note that the expression (\ref{presque-donaldson2}) makes sense since the term $
\left(\sum_{i} |\widetilde{S_{i,k-1}}|^2 \right)^{-1}$ can be considered as a section of $K_M^r \otimes \overline{K_M}^r$. 
\par \label{Remarque-importante} In order to obtain the $\Ric^{-k}(\o_0)^n$-balanced metric $H_{k,r}$ one needs to iterate the operator $Hilb_{\omega_{k,r}}\circ FS$, which gives a sequence $(H_{k,r,p})_{p\in \mathbb{N}}\in Met(H^0(M,L^r))$ initialized with $H_{k,r,0}=H_{k-1,r}$. Here, we denote
$H_{k,r}=H_{k,r,\infty}$.  Remark now that if one expects the algorithm to be convergent, and
thus $H_{k-1,r}$ to be close to $H_{k,r}$ for large $k$ and $r$, then it is natural to assume that $H_{k-1,r,1}$ is close to $H_{k-1,r,\infty}=H_{k,r,0}$, i.e just one step is sufficient to get the $\Ric^{-k}(\o_0)^n$-balanced metric. This justifies at least formally the definition of canonically balanced metrics for Fano manifolds that appeared in \cite[Section 2.2.2]{D3}. 
\begin{definition}
Let $M$ be a Fano manifold and $r \in \mathbb{N}^*$ sufficiently large. For a metric $G\in Met(H^0(M,-K_M^r)$, we define the operator $\widetilde{T}$ on $Met(H^0(M,-K_M^r)$ by
\begin{equation}\label{donaldson}
\widetilde{T}(G)_{ij}= \int_M \frac{\langle S_i , S_j \rangle }{\left(\sum_{i} |S_i|^2\right)^{1+1/r}}  
\end{equation}
for $(S_i)\in H^0(M,-K_M^r)$ a $G$-orthonormal basis. A canonically balanced metric $G_{can}\in Met(H^0(M,-K_M^r))$ is a fixed point of the operator $\widetilde{T}$.
\end{definition}
The operator $\widetilde{T}$ seen as acting on Bergman type metrics is actually a finite dimensional approximation of the $\Ric^{-1}$ operator. Hence, when they exist a priori, we expect the behaviour of the sequence of canonically balanced metrics to be understood via the iterations of the $\Ric^{-1}$ operator, and thus to be related to the behaviour of the K\"ahler-Ricci flow. In that direction, we can show
\begin{theorem}
Assume that $M$ is a Fano Einstein manifold with no nontrivial holomorphic vector field.
Then for $r$ sufficiently large, there exists a sequence of canonically balanced metrics $G_{can,r}\in Met(H^0(M,-K_M^r))$ such that $c_1(FS(G_{can,r}))$ converges when $r$ tends to infinity to the \KE metric $\o_{KE}$.
\end{theorem}
\begin{proof}
The proof is similar to the proof of the main result of  \cite{D1}. Actually, a canonically balanced metric $h_{can,r}\in Met(-K_M^r)$ satisfies  that the Bergman kernel is constant on the manifold, i.e for all $p\in M$,
$$\sum_{i=1}^{N_r} |S_i|^2_{h_{can,r}}(p)= \frac{N_r}{V}$$
where the sections $S_i\in H^0(M,-K_M^r)$ are orthonormal with respect to the inner product
$$<a,b> = \int_M (h_{can,r})^{\frac{r+1}{r}}\otimes a \otimes \overline{b}, $$ and $N_r=\dim H^0(M,-K_M^r)$.  Note that one can generalize the results of Z.\ Lu \cite{Lu} on the asymptotic of the Bergman kernel in that case. The details will appear in a forthcoming paper where the relationship with G.I.T stability will be studied.
\end{proof}
In our implementations we use this notion, i.e we consider the iterations of the operator $\widetilde{T}$. At each step we also choose a normalisation for the metric by fixing essentially its volume. Our tests on toric manifolds have shown that the sequence of metrics defined by  (\ref{presque-donaldson2}) and (\ref{donaldson}) have similar behaviours (i.e iterating many times $Hilb_{\omega_{k,r}}\circ FS$ or just once) and converge to the \KE metric when it exists a priori. These procedures have the advantage to skip the computation of the determinant of the Fubini-Study metric as required by the original notion of balanced metric \cite{D1,Luo,Zha}. Therefore they are clearly much more efficient and we expect that similar methods could be developed for Einstein non Fano manifolds. Finally, remark that even in the case of $\mathbb{CP}^1$ with the anticanonical polarisation, the sequence of balanced and canonically balanced metrics converge at different speeds towards the Fubini-Study metric \cite{D3}. 

\subsection{The case of the projective plane blown up in three points}

Let us consider the toric Fano manifold $\mathbb{P}^2$ blown up in 3 (non aligned) points. From a result of J.\ Song, its $\alpha$-invariant is 1 and thus it possesses a \KE metric (this is also a consequence of Tian's work of classification of Einstein Del Pezzo surfaces). Let us mention that the \KE metric on this manifold has been very recently studied in \cite{DHHKW} by simulating the Ricci flow with PDE techniques.

We implement our algorithms (i.e in order to find balanced and canonically balanced metrics) using the special symmetries of this manifold. The computations of the points on the manifolds are relatively quick since we are essentially reduced to a 2-dimensional real manifold and there are some symmetries given by the action of $\ZZ_2$ (reflections) and $\ZZ_6$ (rotations). The fan of this toric variety is given by the six rays spanned by 
\begin{eqnarray*}
v_0=(1,0), v_1=(1,1), v_2=(0,1), v_3=(-1,0), v_4=(-1,-1), v_5=(0,-1).
\end{eqnarray*}
and as it is well known the polytope is actually the hexagon. There are different ways to choose the points on
this toric manifold but we decided to just generate the points on one of the 6 affine charts
associated to the cone formed by pairs of rays $(v_i,v_{i+1})$ (i.e. by defining a certain cut-off function). \bigskip

\par
\begin{figure}
\begin{center}
\includegraphics[width=5.8cm]{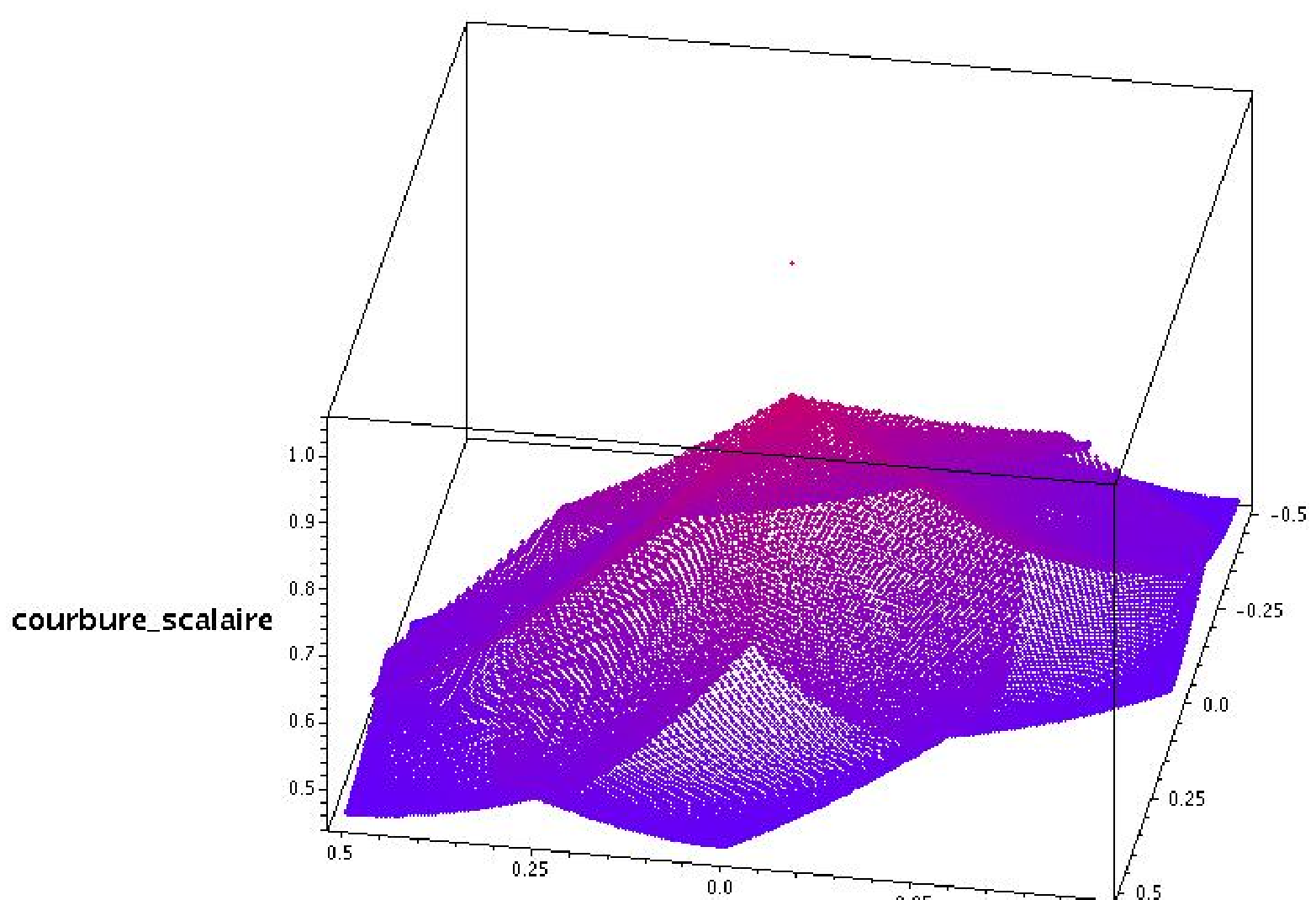}\hspace{0.15cm}
\includegraphics[width=5.8cm]{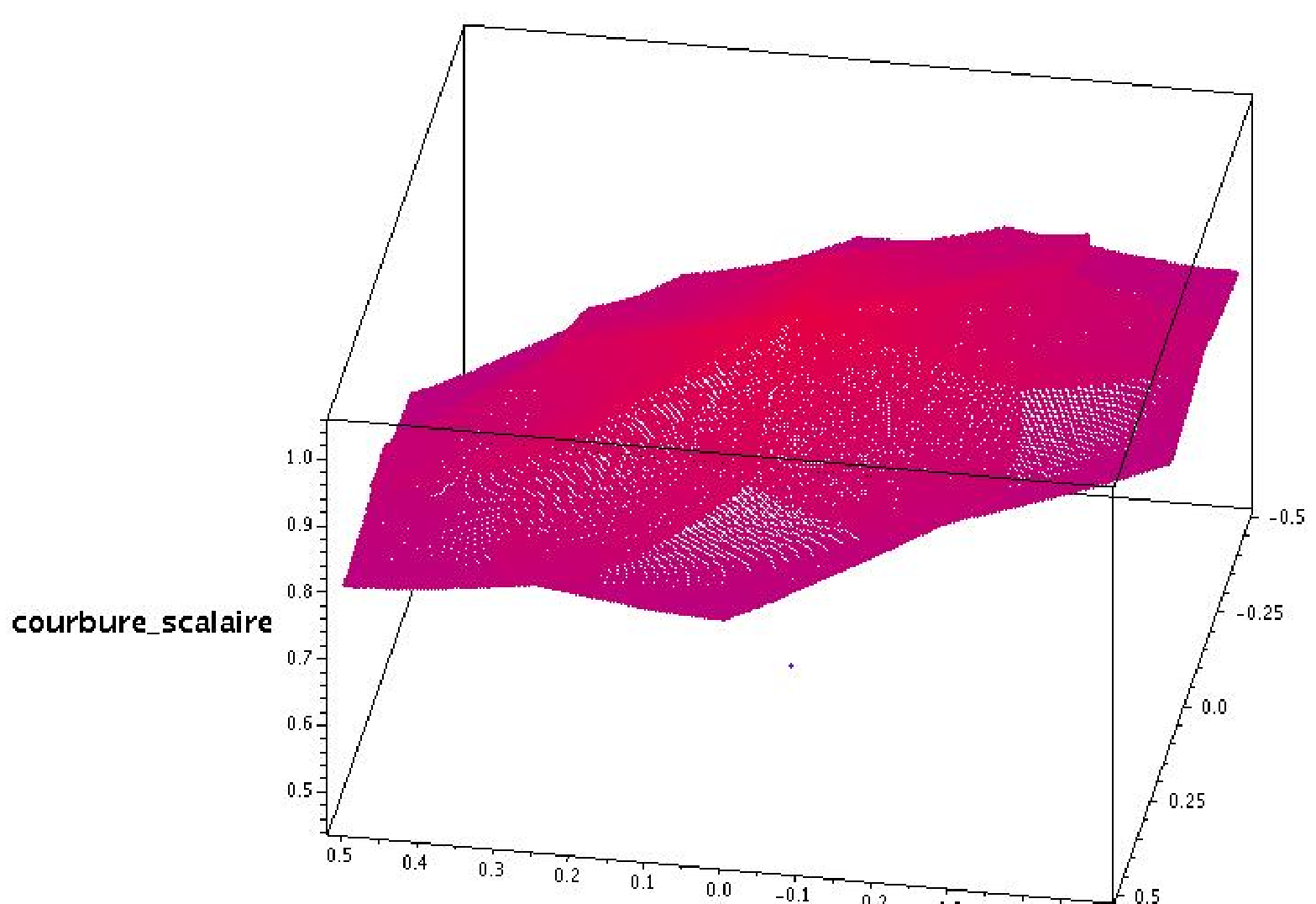}\vspace{0.15cm}\\
\includegraphics[width=5.8cm]{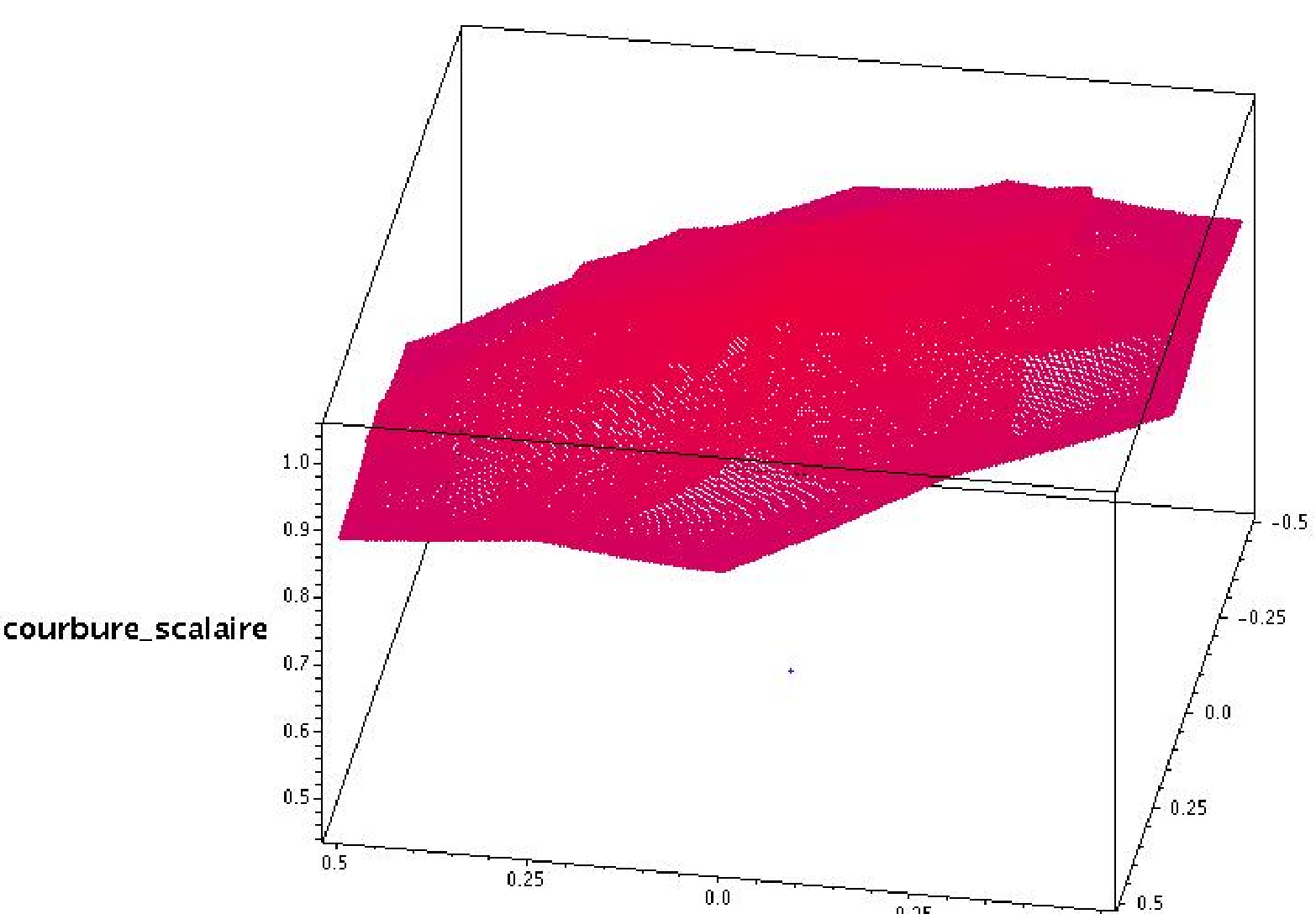}\hspace{0.15cm}
\includegraphics[width=5.8cm]{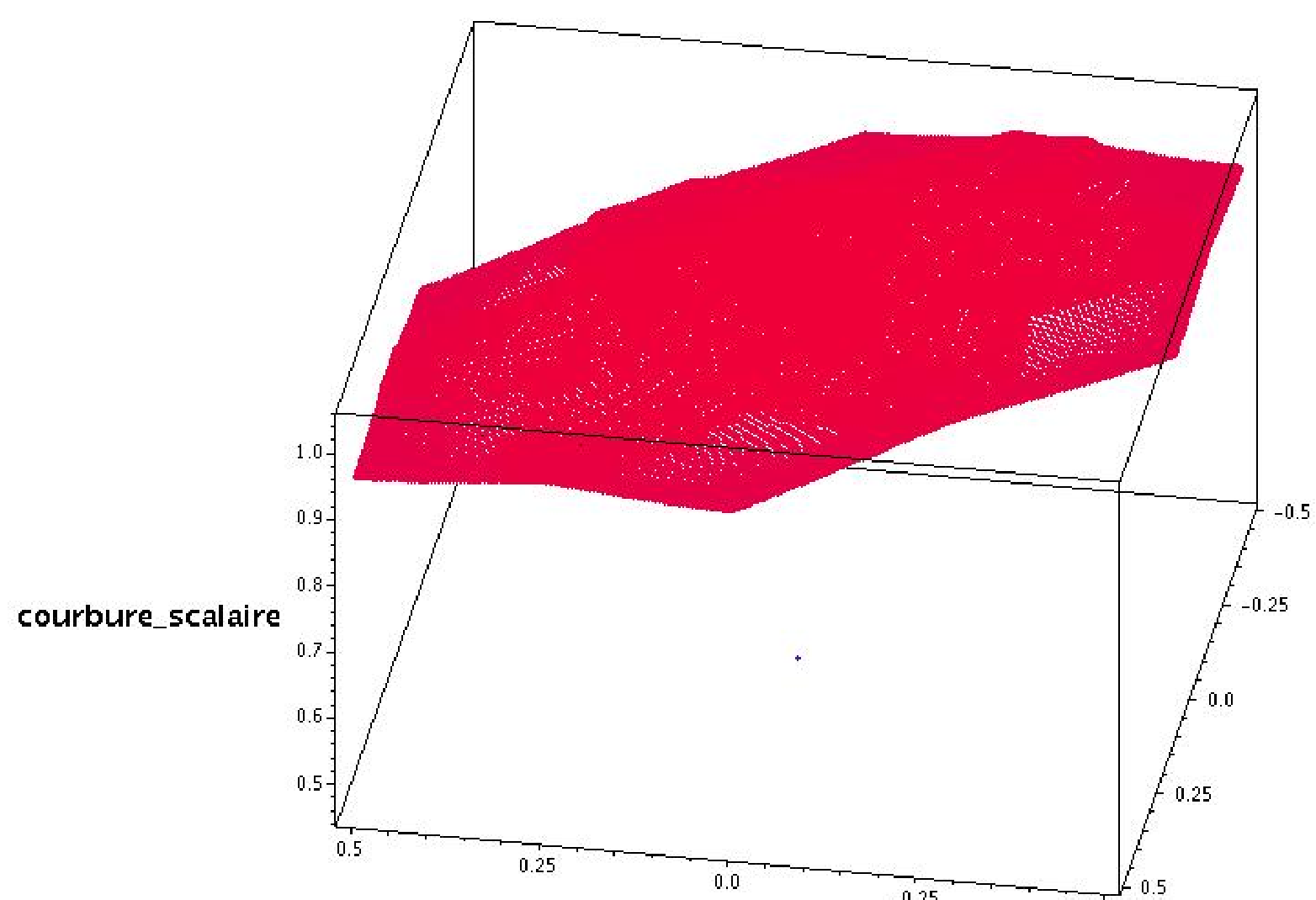}
\end{center}
\caption[$\mathbb{P}^2$ blown up in 3 points]{Iterations no. 10,16,20,40 for the canonically balanced metric}
\label{figures}
\end{figure}

 Our program is written in C++ (compiler gcc 3.4.6) and can be launched essentially with 4 different algorithms. For each algorithm we print the scalar curvature of the new computed metric computed at each point of the manifold. This gives a picture with different colours and one can easily see after some iterations that the scalar curvature varies essentially for the points close to the edges. The computation of the scalar curvature is possible with exact precision (up to the machine precision) since our metrics are algebraic. It has the disavantage to take time since it involves derivatives of order 4, but on the other hand it gives a visual output of our work, see Figure (\ref{figures}). Some animated pictures generated by the program and the program itself can be downloaded from the website of the author\footnote{http://www.ma.ic.ac.uk/$\sim$jkeller/Julien-KELLER-progs.html \\ Some other programs for other Fano Einstein surfaces will be available at this address.}. Despite of this loss of time, all the four algorithms (for the given parameters below) can be run in 1 minute or less on a decent desktop computer. This proves the efficiency of the methods and let us hope that it is possible to compute \KE metric on 3-folds with few symmetries. We now describe the results for each algorithm.
\par For the computation of the balanced metric as defined in \cite{D1,D2,D3}, we choose the parameter $r=8$ (see Section \ref{finite}) and compute approximatively $10^4$ points. After $50$ iterations, the average scalar curvature on the manifold is $0.95$ and the maximum error is $16\%$. 
\par For the computation of the canonically balanced metric as defined in \cite{D3} or our discretization of the Ricci flow (\ref{presque-donaldson2}), we fix again $r=8$ and 
compute approximatively $5\cdot 10^4$ points. After $35$ iterations, the average scalar curvature on the manifold is now $0.99$ and the maximum error is less than $4\%$. 
\par Finally we try to improve the our first two algorithms by using the metric at rank $r$ to compute the metric at rank $r+1$. This is based on a very simple argument that we describe now. If one knows the balanced $h_{r}\in Met(L^r)$ (or canonical balanced metric) at rank $r$, then using the asymptotic of the Bergman function for higher tensor powers (see \cite{Lu} for details) one can write on this surface
$$\frac{N_r}{V}=\sum_{i} |S_i|^2_{h_{r}}(p) = r^2+r\frac{scal(c_1(h_{r}))(p)}{2}+\Gamma(p)+O(1/r)$$
where $\Gamma(p)$ is a certain function which in fact is an algebraic expression of the curvature of $c_1(h_{r})$ and its derivatives. Here $(S_i)_{i=1,..,N_r}$ is an orthonormal basis of $H^0(M,L^r)$ with respect to the $L^2$-metric corresponding to the choice of our algorithm. Once we have computed $h_{r}$, it is clear that we can deduce the value of $\Gamma$ at each point of the manifold.
Now at rank $r+1$, we look for a metric $\tilde{h}_{r+1}$ such that
\begin{eqnarray*}
\sum_{i=1}^{N_{r+1}} |S_i|^2_{\tilde{h}_{r+1}}(p) &=& \frac{N_{r+1}}{V} + \Gamma(p)  \\
&=& \frac{1}{V} \left(N_{r+1}+ \frac{r}{2}\left(1-scal(c_1(h_{r}))(p) \right)\right)
\end{eqnarray*}
with respect to the corresponding $L^2$-metric. Roughly speaking, it corresponds to force the algorithm to get a metric with constant scalar curvature up to an error of size $O\left({1}/{(r+1)^2}\right)$ (instead of only $O\left({1}/{(r+1)}\right)$ for the balanced metric). We apply also the same trick for the canonically balanced metrics. Finally we call these new sequence of metrics as 1-step recursively balanced (1-s.bal in short) or 1-step recursively canonically balanced  (1-s.c.bal in short). \\

The advantage of this method is that it is particularly simple (at least for dimension 2) to program it, since we have already coded the computation of the scalar curvature. On the other hand it asks of course a little bit more time of computation at each step. 
We give now an overview of the results for $r=4$, a choice of $10^4$ points on the manifolds and 15 iterations for each algorithm.\\
\begin{center}
\begin{tabular}{l||c|c|c|c}
 Method $(r=4)$ &                 Balanced & 1-s.bal & Can. balanced   & 1-s.c.bal \\
  Avg scalar curvature &  0.786    & 0.888   &   0.949    &  0.984\\
 Max scalar curvature &   1.079    & 1.041   &   1.012    &  1.041\\
 Min scalar curvature &   0.618    & 0.717   &   0.820    &  0.867\\
  Time (sec)           &     9.1     &  16.1   &     8.0    &   15.5  \\
\end{tabular}
\end{center}

\noindent As an example, we obtain Figure \ref{figures} p.\pageref{figures} for the canonically balanced metric with $r=12$ and various iterations. One can see that for $40$ iterations we get a metric with scalar curvature almost equal to $1$ everywhere. Using this metric, one can find numerical approximations\footnote{http://www.ma.ic.ac.uk/$\sim$jkeller/Julien-KELLER-progs.html} of geodesics for the \KE  metric. It seems to give a numerical evidence that the geodesic equations on this manifold form an integrable system.

\section{Iteration for other classes}

\subsection{The case of negative first Chern class}
Assume now that $K_M>0$ i.e that $M$ has negative first Chern class. It is well known by \cite{Ya} that such manifolds have a unique \KE metric $\oKE$ such that
$$\Ric(\oKE)=-\oKE.$$
 Let $\Ric\o+\o=\i\ddbar h$ with $\V\intM e^{h}\o^n=1$.
In that case, the naive discretization of the normalised K\"ahler-Ricci flow leads us to define the iteration scheme by
\begin{equation}\label{negfirstchernclass}
-\Ric\o_l=2\o_l -\o_{l-1}.
\end{equation}
This can be rewritten in terms of the following \MA equations
\begin{equation}
\ovpln=e^{2\varphi_l-\varphi_{l-1}}\o_{\sum_{j=1}^{l-1}\varphi_j}^n,\qquad \V\intM e^{2\varphi_l-\varphi_{l-1}}\o_{l-1}^n=1,\q \all\,l\in\NN.\label{ConenegEq2}
\end{equation}
and we will denote as previously $\Phi_l=\sum_{j=1}^l \varphi_l$ and $\o_l=\o+\sqrt{-1}\ddbar \Phi_l$. Note that by the general results on \MA equations \cite{Ko}, for any $l$, there exists a smooth solution $\varphi_l$ of (\ref{ConenegEq2}) and thus the iteration scheme is well-defined. By concavity of the $\log$ applied to (\ref{ConenegEq2}), we have directly 
\begin{lemma}\label{decrease_c1_neg}
Let $\mu=-1$. Then 
\begin{eqnarray}
E_0(\o_{l-1},\o_l)= & -J(\o_l,\o_{l+1})+\V\intM (2\varphi_l-\varphi_{l-1}) \o_{l-1}^n 
\leq 0&
\end{eqnarray}
i.e the iterations (\ref{ConenegEq2}) decrease the functional $E_0$.
\end{lemma}

\begin{proposition}\label{c0estimate}
Let $p\in(1,\infty)$. There exists a constant $C_1$ depending only on $p,M$ and $\o$ such that
$$
\n{\Phi_l}\LinfMo\leq C_1,\q \all\,l\in\NN.
$$
\end{proposition}

\begin{proof}
Let $p\in(1,\infty)$. According to the work of \cite{Ko} 
(or \cite{Bl} if $p\in(2,\infty)$ which will suffice for the proof) and in view of (\ref{ConenegEq2}) it suffices to prove that
$$
\n{e^{h+2\Phil-\Phi_{l-1}}}\LpMo\leq C_2,\q \all\,l\in\NN,
$$
for $C_2=C(M,\o,p)$. From the normalisation in (\ref{ConenegEq2}) it follows that
$\V\intM (h+2\Phil-\Phi_{l-1}) \on\le0$ where $h$ denotes the Ricci deviation of $\o$. In particular 
$\sup_M (h+2\Phil-\Phi_{l-1})\leq -\VintM (h+2\Phil-\Phi_{l-1}-\sup_M(h+2\Phil-\Phi_{l-1}))\on$. Therefore
\begin{eqnarray*}
\VintM e^{p(h+2\Phil-\Phi_{l-1})}\on
&\le&
e^{p\sup_M(h+2\Phil-\Phi_{l-1})} \\
&\le&
e^{-p\VintM (h+2\Phil-\Phi_{l-1}-\sup_M(h+2\Phil-\Phi_{l-1}))\on} \\
&\le&
e^{p\, \osc h}e^{-p\VintM (2\Phil-\Phi_{l-1}-\sup_M(2\Phil-\Phi_{l-1}))\on}.
\end{eqnarray*}
Let $G$ be the Green function for $\D=\D_{\dbar}$ satisfying $\intM G(x,y)\on(y)=0$.
Since $\o_0+2\Phil-\Phi_{l-1}=2\o_l-\o_{l-1}>0$ the Green formula gives
\begin{eqnarray*}
\lefteqn{\VintM (\sup_M(2\Phil-\Phi_{l-1})-(2\Phil-\Phi_{l-1}))\on}\\&=&-\VintM G(x_0,y)\D(2\Phil-\Phi_{l-1})\on(y) \leq nA
\end{eqnarray*}
where $A=-\inf_{M\times M} G(x,y)$ and $x_0\in M$ satisfies $(2\Phil-\Phi_{l-1})(x_0)=\sup_M(2\Phil-\Phi_{l-1})$. Therefore $C_2=e^{p(\osc h+nA)}.$
\end{proof}

Combining the previous proposition and Lemma \ref{decrease_c1_neg}, we obtain
\begin{corollary}
Let $M$ be a manifold with $K_M>0$. Then, the sequence of K\"ahler  metrics defined by (\ref{negfirstchernclass}) converges to the \KE metric $\oKE$ on $M$.
\end{corollary}

Actually, another iteration scheme is also very natural, especially if we assume the metrics $\o_l$ and $\o_{l+1}$ from  (\ref{negfirstchernclass}) to be very \lq close\rq\, for a large $l$. 
We can define the iteration scheme 
\begin{equation}\label{negfirstchernclass2}
-\Ric\widetilde{\o}_l=\widetilde{\o}_{l-1}.
\end{equation}
which means that we solve at each step the following \MA equations
\begin{equation}
\widetilde{\o}_{\sum_{j=1}^l \widetilde{\varphi}_j}^n=e^{h+\sum_{j=1}^{l-1} \widetilde{\varphi}_j}\on,\qquad \V\intM e^{h+\sum_{j=1}^{l-1} \widetilde{\varphi}_j}\on=1,\q \all\,l\in\NN.\label{ConenegEq}
\end{equation}
 Let's call $\widetilde{\Phi}_l=\sum_{j=1}^{l} \widetilde{\varphi}_j$. Up to our knowledge,  it is not clear whether the iteration scheme (\ref{negfirstchernclass2}) decreases any functional. Now, one can check that the proof of Proposition \ref{c0estimate} can be adapted to that case without must change, i.e one has straightforward a $C^0$-estimate of the potentials $\widetilde{\Phi}_l$ and thus convergence up to extraction of a subsequence. On another hand, the iteration scheme (\ref{negfirstchernclass2}) is more natural from the finite dimensional perspective as we shall explain now. 
\par Similarly to the notion of canonically balanced metric introduced by (\ref{donaldson}) in Section \ref{finite}, one can define a notion of canonically balanced for manifolds with negative first Chern class. Indeed, it turns out to consider the iteration procedure defined by the iterations of the map $\widetilde{T}=\widetilde{Hilb} \circ FS$ where 
$$\widetilde{Hilb}(h)(S_i,S_j)=\int_M h^{\frac{k-1}{k}}\otimes S_i \otimes \overline{S_j}$$
for $h\in Met(K_M^k)$ and $S_i\in H^0(M,K_M^{k})$. Here we see $h^{\frac{k-1}{k}}$ as an element of $K_M^{1-k}\otimes \overline{K_M}^{1-k}$. From our previous discussion p.\pageref{Remarque-importante} Section \ref{Remarque-importante}, one can expect naturally the iterates of $\widehat{T}$ to approximate in the space $H^0(M,K_M^k)$ the discretization of the normalised K\"ahler-Ricci flow given by (\ref{negfirstchernclass2}). \\ 
As explained in \cite{Ts}, another remarkable fact is that one can define an iterative procedure in higher dimensional spaces, i.e by defining now the maps $\widehat{T}=\widehat{Hilb} \circ FS:Met(K_M^k)\rightarrow Met(K_M^{k+1})$ where
$$\widehat{Hilb}(h)(S_i,S_j)=\int_M h\otimes S_i \otimes \overline{S_j}.$$
Let's define $h_k=\widehat{T}^{(k)}(h_0)$ for a metric $h_0\in Met(K_M^{k_0})$ for $k_0$ sufficiently large. As observed by Weinkove and Song \cite[Theorem 1]{SW2}, one has a uniform convergence for $k \rightarrow +\infty$ of the metrics $h_k^{\frac{1}{k+k_0}}$ to the smooth hermitian metric $h_{KE}\in Met(K_M)$ such that $c_1(h_{KE})=\oKE$ is the \KE metric on $M$. Nevertheless this method does not seem to have numerical applications straightforward since it involves to consider higher rank matrices (and thus their inverses). 

\begin{remark}
One can also derive an iterative scheme for manifolds with trivial canonical bundle. This leads to a refinement of the notion of $\nu$-balanced metric \cite{D3}, and is useful for our numerical study of Ricci-flat metrics on Calabi-Yau threefolds. This will be addressed in a forthcoming paper.
\end{remark}

\subsection{The case of non-canonical classes}

 Let $\o\in\Ho$ be a \K representative of an arbitrary
class $\O$ in the \K cone of a Fano manifold $M$. Let $\sigma\in\hc0$ be a representative
of $c_1(M)$. By the Calabi-Yau theorem there exists a unique \K representative $\o_1$ of
${\O}$ whose Ricci form equals $\sigma$. Define a map $${\RicmO}: \hc0 \ra \Ho$$ by
$\RicmO\sigma:=\o_1$ (when appropriate
norms are chosen, it defines again a homeomorphism of Banach spaces).
We may therefore define an operator $\Ric_{[\o]}$ on $\Ho$
by $$\Ric_{[{\o}]}:=\RicmO \circ\Ric\circ\Ric.$$ This gives rise to an iteration on $\O$ by setting
$\rico k:=\RicmO\circ\ric k\circ\Ric$ for each $k\in\ZZ$. Note that this induces
again a filtration on $\Ho$ defined by $\ho k:=\{\o\in\Ho:
\Ric\o\in\hc k\}$.

From what we have proved previously,  we expect that the orbits of the new dynamical systems converge, if and only
$M$ is \KEno, to a \K representative which is characterised by the property that its
Ricci form is \KEno. For $\O\ne c_1(M)$, one can ask how is this dynamical system related
to the study of the space $\Ho$. In particular, it would be interesting to relate this to the
existence problem of extremal metrics for \K classes near $c_1(M)$.\\
\vspace{1cm}
\\
\noindent {\small {\bf Acknowledgements.} The author is extremely grateful to S.K. \ Donaldson for enlightening conversations and his constant encouragement. He also thanks R.\ Bunch for his very useful help with C++ programs. This work has also benefited from stimulating discussions with X.X.\ Chen, D.\ Panov, D.H.\ Phong, S.\ Simanca, G.\ Szekelyhidi and R.P.\ Thomas. Finally he thanks P. Eyssidieux for his support throughout the years.}

\textsc{\\
{\bf Julien Keller}\\
Imperial College, London\\ 
j.keller@imperial.ac.uk
}
\end{document}